\documentclass{article}

\usepackage{arxiv}

\usepackage{amsmath}
\usepackage{hyperref}
\usepackage{algorithm}
\usepackage{xcolor}

\usepackage{lipsum}
\usepackage{amsfonts}
\usepackage{graphicx}
\usepackage{epstopdf}
\usepackage{algpseudocode}
\usepackage{multirow}

\usepackage{amssymb}
\usepackage{bm}

\newcommand{\bq}{{\bm q}}
\newcommand{\bs}{{\bm s}}
\newcommand{\bv}{{\bm v}}

\newcommand{\bx}{{\bm x}}
\newcommand{\by}{{\bm y}}

\newcommand{\bbR}{\mathbb R}

\newcommand{\calB}{\mathcal{B}}
\newcommand{\calD}{\mathcal{D}}

\newcommand{\calS}{\mathcal{S}}

\newcommand{\raw}{\rightarrow}

\newcommand{\vn}[1]{\left|\left|#1\right|\right|}
\newcommand{\msdi}[1]{\texttt{MSD}_{\,i}\left(#1\right)}

\newcommand{\as}{\text{$\langle\texttt{samp}\rangle$}}
\newcommand{\qpdf}{\texttt{qpdf}}

\newcommand{\setalglineno}[1]{%
  \setcounter{ALG@line}{\numexpr#1-1}}

\ifpdf
  \DeclareGraphicsExtensions{.eps,.pdf,.png,.jpg}
\else
  \DeclareGraphicsExtensions{.eps}
\fi

\usepackage{enumitem}
\setlist[enumerate]{leftmargin=.5in}
\setlist[itemize]{leftmargin=.5in}

\title{Data-driven geometric scale detection via Delaunay interpolation\thanks{Submitted to the editors March 24 2022.
\funding{This work was funded by LLNL--LDRD Project No.\ 21--ERD--028}}}

\author{Andrew Gillette\thanks{Lawrence Livermore National Laboratory, Livermore, CA  
  (\email{gillette7@llnl.gov}, \url{https://akgillette.github.io/}).}
\and Eugene Kur\thanks{Lawrence Livermore National Laboratory, Livermore, CA 
  (\email{kur1@llnl.gov}).}}

\usepackage{amsopn}

\usepackage[utf8]{inputenc} 
\usepackage[T1]{fontenc}    
\usepackage{url}            
\usepackage{booktabs}       
\usepackage{amsfonts}       
\usepackage{nicefrac}       
\usepackage{microtype}      
\usepackage[capitalize,nameinlink]{cleveref}       
\usepackage{lipsum}         
\usepackage{graphicx}
\usepackage{doi}

\title{The Delaunay Density Diagnostic}

\date{}

\author{ \href{https://orcid.org/0000-0002-1195-5924}{\includegraphics[scale=0.06]{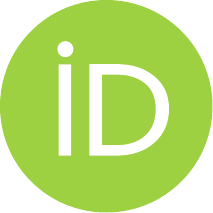}\hspace{1mm}Andrew Gillette}\thanks{\url{https://people.llnl.gov/gillette7}} \\
	Center for Applied Scientific Computing \\
	Lawrence Livermore National Laboratory \\
	7000 East Avenue, Livermore, CA 94550 \\
	\texttt{gillette7@llnl.gov} \\
	\And
	\href{https://orcid.org/0000-0002-1831-3887}{\includegraphics[scale=0.06]{orcid.pdf}\hspace{1mm}Eugene Kur} \\
	Strategic Deterrence \\
	Lawrence Livermore National Laboratory \\
	7000 East Avenue, Livermore, CA 94550 \\
	\texttt{kur1@llnl.gov} \\
}

\renewcommand{\headeright}{}
\renewcommand{\undertitle}{}
\renewcommand{\shorttitle}{Data-driven geometric scale detection via Delaunay interpolation}

\hypersetup{
pdftitle={Data-driven geometric scale detection via Delaunay interpolation},
pdfsubject={mathematics, computer science, statistics},
pdfauthor={A. Gillette and E. Kur},
pdfkeywords={function approximation, scale detection, interpolation, Delaunay triangulation},
}

\newenvironment{@abssec}[1]{%
     \if@twocolumn
       \section*{#1}%
     \else
       \vspace{.05in}\footnotesize
       \parindent .2in
         {\upshape\bfseries #1. }\ignorespaces 
     \fi}
     {\if@twocolumn\else\par\vspace{.1in}\fi}
     
\newcommand\AMSname{AMS subject classifications}
\newenvironment{AMS}{\begin{@abssec}{\AMSname}}{\end{@abssec}}

\newcommand\keywordsname{Key words}
\renewenvironment{keywords}{\begin{@abssec}{\keywordsname}}{\end{@abssec}}

\Crefname{figure}{Figure}{Figures}

\newtheorem{claim}{Claim}

\begin{document}
\maketitle



\begin{abstract}
Accurate approximation of a real-valued function depends on two aspects of the available data: the density of inputs within the domain of interest and the variation of the outputs over that domain.
There are few methods for assessing whether the density of inputs is \textit{sufficient} to identify the relevant variations in outputs---i.e., the ``geometric scale'' of the function---despite the fact that sampling density is closely tied to the success or failure of an approximation method.
In this paper, we introduce a general purpose, computational approach to detecting the geometric scale of real-valued functions over a fixed domain using a deterministic interpolation technique from computational geometry. The algorithm is intended to work on scalar data in moderate dimensions (2-10).
Our algorithm is based on the observation that a sequence of piecewise linear interpolants will converge to a continuous function at a quadratic rate (in $L^2$ norm) if and only if the data are sampled densely enough to distinguish the feature from noise (assuming sufficiently regular sampling).
We present numerical experiments demonstrating how our method can identify feature scale, estimate uncertainty in feature scale, and assess the sampling density for fixed (i.e.~static) datasets of input--output pairs.
We include analytical results in support of our numerical findings and have released lightweight code that can be adapted for use in a variety of data science settings.


\end{abstract}



\section{Introduction}
\label{sec:intro}



Identification of sufficient sampling density is an essential and ongoing challenge in data science and function modeling.  
For any problem context, too little data raises concerns of over-fitting while too much data risks under-fitting and inefficient computational pipelines.
While theorems and error estimates can provide rough bounds on requisite sampling density, more often density is selected by heuristics, trial and error, or rules of thumb.

In this work, we introduce the \textit{Delaunay density diagnostic}: a computational technique that can help identify the presence and scale of geometric features in a function $f:\bbR^d\raw\bbR^1$ as values of $f$ are iteratively collected.
By geometric features we mean any non-linear behavior that is not noise, and hence the presence of a feature requires sufficient sampling density to distinguish it from noise.
More formally, the goal of the diagnostic is to provide a robust computation of the rate at which a piecewise linear interpolant of $f$ changes as additional sample points are incorporated in batches of user-specified size.

A simple illustration of the challenges present in this effort can be seen by examining the Griewank function on $\bbR^1$ ($d=1$), defined by
\[g_1(x) := \frac{1}{4000}x^2 -\cos(x) + 1.\]
At a large scale, e.g.~$x\in[-10^4, 10^4]$, the quadratic term is dominant and only an extremely dense sampling of the interval would be able to distinguish the cosine term from random noise.
At a small scale, e.g.~$x\in[-1,1]$, the cosine term is dominant and even a very dense sampling of the interval would be insufficient to pick up the quadratic term behavior.
At intermediate scales, the interplay between sampling density and feature representation is more subtle, motivating the need for a computational approach to the problem.

We first describe the Delaunay density diagnostic in the simplest context and in a manner slightly distinct from how it is actually implemented.  
Assume we can rapidly compute the output of a function $f:\bbR^d\raw\bbR^1$ for any $\bx\in\bbR^d$.
For now, we assume there is no uncertainty or significant numerical error in the output computation, as can be the case when $f$ is given by an analytical formula or a trained neural network of mild complexity.
The following must then be specified by the user:
\begin{itemize}
	\item A $d$-dimensional \textbf{bounding box} $\calB\subset\bbR^d$. All data and computations will be confined to this box. In practice, it is best to rescale the box to a standard hypercube (such as $[0,1]^d$), which would induce a similar rescaling on any data contained within.
	\item An initial set of \textbf{sample points} $\{\bx_j\}\subset\calB$. Ideally, these points are spread uniformly over $\calB$. The initial Delaunay interpolant will be determined by this set. 
	\item An \textbf{upsampling growth factor} $b\in (1,2]$.  Sample points will increase by $b^d$ in each iteration, typically via uniform upsampling.  Each iteration determines a finer-grained Delaunay interpolant than the previous.
	\item A set of \textbf{query points} $\{\bq_i\}\subset\calB$. The Delaunay interpolant will be evaluated at these points. To avoid too many query points falling outside the convex hull of the sample points (thus resulting in geometric extrapolation), the query points should be supported in a proper subset of $\calB$. In general, the larger $d$ is, the smaller the subset bounding the query points must be, but the nature of the sample points also plays a role (e.g.~ if the samples are intentionally clustered in some region).
\end{itemize}
In the simplest case, we choose the query points to be a regular, axis-aligned lattice of $p^d$ points forming a cube in $\bbR^d$, for some $p$, and $\calB$ to be a box with the same center as the lattice but larger diameter.
The main computation is an iterative procedure:
\begin{enumerate}
\item For each query point $\bq_i$, compute and store ${\hat f}(\bq_i)$, the unique Delaunay piecewise-linear interpolant defined by $\{\bx_j,f(\bx_j)\}$, evaluated at $\bq_i$.
\item Let $\|\{\bx_j\}\|=n$.
Let $\ell$ be the integer closest to $\left[ b~n^{1/d}- (b-1)\right]^d -n$.
Randomly sample $\ell$ more points from $\calB$ and add these points to the set of samples so that $\|\{\bx_j\}\|=n+\ell$.
\item Repeat until $\|\{\bx_j\}\|$ exceeds a specified threshold.
\end{enumerate} 
If some $\bq_i$ lies outside the convex hull of the $\{\bx_j\}$ at a step of the process, the Delaunay interpolant is not defined and we just store a \texttt{nan} value. 
This case can be avoided by increasing the diameter of $\calB$ or increasing the initial number of points in $\{\bx_j\}$.

Observe that the values of ${\hat f}(\bq_i)$ may change at each step, since the introduction of new samples may change the sample points defining the computation of ${\hat f}(\bq_i)$. 
We use this observation to compute an approximate rate of convergence of the interpolants to the true function $f$ over the region defined by the query points. 
Let ${\hat f}_k(\bq_i)$ denote the interpolated value at $\bq_i$ at the $k$th step of the iterative procedure.
For each $k\geq 2$, we compute a rate $r_k$ defined by 
\[r_k := \log_{\,b} \left(\frac{\sqrt{\msdi{{\hat f}_{k-1}(\bq_i) - {\hat f}_{k-2}(\bq_i)}}}{\sqrt{\msdi{{\hat f}_{k}(\bq_i) - {\hat f}_{k-1}(\bq_i)}}}\right).\] 
Here, $\msdi{\cdot}=\frac{1}{|\{\bq_i\}|}\sum_i(\cdot)^2$ denotes the mean squared difference over the index $i$, which provides a measure of the difference between successive interpolants in the procedure; the square root of the mean squared difference is a discretization of the $L^2$ norm.
The definition of $r_k$ is inspired by the computation of convergence rates for finite element methods over a sequence of increasingly refined meshes~\cite{BS1994}.
We now make the following claims, which will be supported by numerical and analytical results later in the paper:

\begin{claim}
\label{claim:rate2}
If the set $\{(\bx_j,f(\bx_j)\}$ used to define $\hat f_k$ contains sufficient information to reconstruct the geometric features present in $f$ over the set of query points $\{\bq_i\}$ then the mean squared difference rate $r_k$ will be approximately $2$.
\end{claim}

\begin{claim}
\label{claim:rate0}
If the set $\{(\bx_j,f(\bx_j)\}$ used to define $\hat f_k$ cannot distinguish geometric features from random noise over the set of query points $\{\bq_i\}$ then the mean squared difference rate $r_k$ will be approximately $0$.
\end{claim}

We envision the results of these claims being relevant to the multitude of scientific machine learning problems in contemporary literature where the goal is to approximate some unknown function based on unstructured numerical data.  
Many such problems have inputs in $\bbR^d$ for $d>2$, making visualization difficult, but not $d\gg 2$, making techniques designed for high-dimensional data not necessarily applicable.  
In such settings, there are few computational techniques for robust identification of sufficient sampling density and hence the proposed diagnostic (which works best for $2<d<10$) could aid in determining whether a given sampling procedure is insufficient, sufficient, or excessive in the context of a specific problem or application.
We note, however, that sampling density is important throughout computational mathematics and the technique has no intrinsic relation to machine learning.

The diagnostic effectively evaluates nonlinear feature presence at different scales, but has several limitations that need to be considered:

\begin{itemize}
	\item The Delaunay triangulation becomes expensive to construct in higher dimensions, even partially, which limits the maximum dimension to around 10.
	\item Much of our work on the diagnostic assumes (nearly) uniform sampling for the sample points. For non-uniform sampling, the ``features'' detected by the diagnostic would be a mix of features of the function and the sampling density. This could potentially be resolved by using a norm adjusted for the sampling density, but we leave such considerations for future work. If the sampling density is very non-uniform (such as if data lies on a lower-dimensional manifold), this could significantly degrade the performance of the Delaunay interpolant and prevent computation of the diagnostic.  Of note, the software we will use for computing the Delaunay interpolant can sometimes detect if the data lies on or nearly on a lower-dimensional manifold.
	\item The diagnostic does not distinguish between geometric features present in all coordinates and those present in only some coordinates. For example, suppose a 4-dimensional function only depended on two of the four coordinates and the features had characteristic size $\ell$. For a bounding box with side length $L$, the diagnostic would require somewhat more than $(L/\ell)^4$ sample points to resolve the features, while an algorithm that could take advantage of the lower dimensionality could fit the function well with only somewhat more than $(L/\ell)^2$ points. In such cases, the diagnostic is too conservative, claiming more points are needed to adequately resolve features even when appropriate algorithms can fully resolve them without additional sampling. 
\end{itemize}



\section{Context in the literature} 
\label{sec:bkgd}

Techniques for fitting a function to unstructured numerical data have been studied in multiple disciplines for many decades. 
Some prominent examples include the Kriging interpolation method~\cite{SWMW1989,SMKM2001,VK2004}, 
radial basis functions~\cite{B2000acta,B2003},
response surface methodology~\cite{KM2010}, and 
generalized linear models~\cite{NW1972,MN2019}.
Interpolation based on ``nearest neighbors'' is appealing for its simplicity, but has not emerged as a feasible competitor to the above methods.
In general, methods have been developed to serve specific needs in engineering and statistics communities but have trouble scaling to the size and dimension of modern datasets.


Perhaps the most closely related research area---and in fact the inspiration for this work---is the notion of data oscillation in finite element methods (FEM); see e.g.~\cite{MNS2000}.
In FEM, the goal is to approximate the solution $u$ to a partial differential equation $Du = f$ by creating a piecewise polynomial with respect to a mesh of the domain.
Data oscillation, $osc_{\,\mathcal{M}}$, is a quantity that measures the variation of $f$ over a fixed mesh $\mathcal{M}$.
If $osc_{\,\mathcal{M}}$ is large, $f$ has fine scale geometric structures that cannot be resolved by approximation with respect to $\mathcal{M}$, meaning a finer mesh must be used.
Hence, attaining optimal rates of convergence in adaptive FEM (in which meshes are partially refined in an iterative process) depends on controlled data oscillation; control can be attained  by assumption, by luck, or by computational detection and mesh refinement, often based on the classical theory of Richardson extrapolation~\cite{R1911,RG1927}.

The issue of function-dependent sampling density requirements is by no means restricted to the finite element world.
Metrics to assess data variation with respect to sample points have been devised for application-specific contexts but have no standard nomenclature.  
These include the ``grid convergence index''\cite{CK1997} and ``index of resolution quality''~\cite{CCY2005} for computational fluid dynamics, ``local feature size'' for homeomorphic surface reconstruction from point clouds~\cite{ACDL2000}, and the ``coefficient of variation'' for LIDAR sampling~\cite{GLYA2010}.
In machine learning communities, various notions of ``holdout sets'' and ``cross validation'' are in vogue~\cite{BKL1999,BG2003,BHT2023}. 
The list could go on. 



This work introduces a methodology for assessing data variation on unstructured numerical data with input dimension $d$ up to $\approx~10$.
Our approach is not meant to generalize or replace the application-specific techniques and metrics described above.
Rather, we aim to provide users working with high volume but relatively low-dimensional numerical datasets a means to assess data variation and sampling density in a rigorous framework, subject to clearly stated assumptions.
Such datasets are now ubiquitous in scientific disciplines, but analyzing function variation with respect to a mesh is often stymied by the so-called ``curse of dimensionality.'' 
For $d>3$, it quickly becomes infeasible to compute, store, or manipulate the complete mesh structure of a collection of unstructured data points.

As we will demonstrate, the lack of scalability of mesh management can be circumvented for interpolation tasks if the interpolated value can be determined using only a sparse subset of an implied---but not computed---mesh data structure.
Delaunay theory provides the requisite mathematical results for an implied mesh structure and the recently developed algorithm \texttt{DelaunaySparse} provides a practical tool for such computations.

\subsection{Delaunay interpolation} 

\begin{figure}[h]
\begin{center}
\includegraphics[width=.3\textwidth]{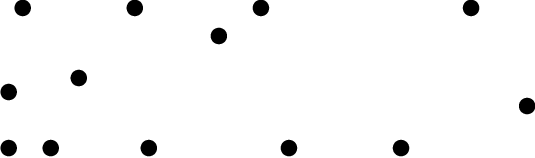}
\qquad\qquad
\includegraphics[width=.31\textwidth]{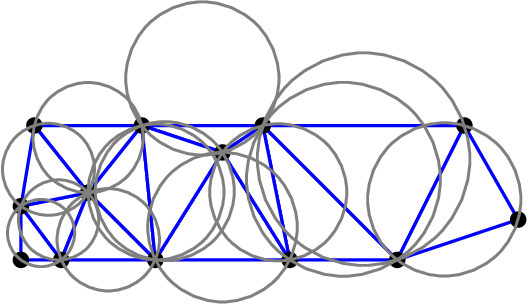} 
\end{center}
\caption{An arbitrary collection of points $\calD$ in $\bbR^2$ (left) has a unique triangular mesh (right, blue lines), called the Delaunay triangulation. 
Every triangle in the Delaunay mesh satisfies the ``empty ball criterion'': the open circumball whose boundary passes through the vertices of the triangle does not contain any points from $\calD$.
The boundaries of the circumballs for the Delaunay triangulation in the figure above are shown as grey circles.
For $d > 2$, these properties generalize to meshes of $d$-simplices and associated $d$-dimensional circumballs.
}
\label{fig:del-intro}
\end{figure}

Let $f: \bbR^d \rightarrow \bbR^\ell$ be a multivariate function whose outputs are known at a collection of $n$ data points $\calD \subset \bbR^d$.
Assume that $\calD$ is truly a $d$-dimensional sample in the sense that it does not lie entirely in a hyperplane of $\bbR^d$; we remark that if there was such a hyperplane, the interpolation algorithm would detect this.
Then, the \textit{convex hull} of $\mathcal{D}$, denoted $CH(\calD)$, is a $d$-dimensional, flat-faced, convex region in $\bbR^d$.
The \textit{Delaunay triangulation}, denoted $DT(\calD)$, is an unstructured mesh of $CH(\calD)$ consisting of $d$-simplices with vertices in $\calD$ that satisfy the open circumball property, which is described in the caption of Figure \ref{fig:del-intro}.
The Delaunay triangulation exists and is unique, except for some corner cases, which can still be handled robustly.
For instance, data stored on a rectilinear lattice does not have a unique Delaunay triangulation, but in such situations, the geometry of the data is known \textit{a priori}.
For a randomly generated data set in $\bbR^d$, as is employed here, the probability of a unique Delaunay triangulation is 1.

The Delaunay triangulation can be used to define a unique piecewise linear interpolant, called the \textit{Delaunay interpolant}, of the values $f(\calD)$.
The Delaunay triangulation is widely considered to be an optimal simplicial mesh for the purposes of multivariate function interpolation \cite{chen2004optimal,rajan1994optimality}, and the Delaunay interpolant (defined below) has been studied in the context of both regression \cite{liu2019nonparametric,lux2019interpolation} and classification \cite{belkin2018overfitting} problems.

Given any point $\bq\in CH(\calD)$, the \textit{Delaunay interpolant} at $\bq$ is defined as follows. 
Let ${\calS}$ be the simplex in $DT(\calD)$ with vertices $\bs_1$, $\ldots$, $\bs_{d+1}\in\calD$ such that $\bq \in {\calS}$.
Then there exist non-negative weights $w_1$, $\ldots$, $w_{d+1}\in\bbR$ such that $\bq = \sum_{i=1}^{d+1} w_i \bs_i$ and $\sum_{i=1}^{d+1} w_i = 1$.
The value of the Delaunay interpolant ${\hat f}_{DT}$ at $\bq$ is given by
\begin{equation}
{\hat f}_{DT}(\bq) := w_1f(\bs_1) + w_2f(\bs_2) + \ldots + w_{d+1}f(\bs_{d+1}).
\label{eq:fDTx}
\end{equation}
The approximation defined by ${\hat f}_{DT}(\bq)$ is a continuous, piecewise linear interpolant of $f$ for $\calD$.
In particular, if $\bq$ lies at the interface of multiple $d$-simplices in $DT(\calD)$, the value of ${\hat f}_{DT}(\bq)$ is not dependent on the choice of simplex used to compute it.

\subsection{The \texttt{DelaunaySparse} algorithm}
\label{subsec:DS}
While computing a data structure for the full $DT(\calD)$ is not computationally feasible beyond very low dimensions, the value of ${\hat f}_{DT}(\bq)$ only requires detection of a simplex in $DT(\calD)$ that contains $\bq$.
\texttt{DelaunaySparse} is a recently developed algorithm and software package that exploits this observation to provide efficient computation of the Delaunay interpolant in high dimensions at a user-provided set of input points.
\cref{alg:DS} outlines the general strategy of \texttt{DelaunaySparse}; further details can be found in \cite{Cetal2020,chang2018polynomial}.


\begin{algorithm}
\caption{\texttt{DelaunaySparse}~\cite{Cetal2020}}
\label{alg:DS}
\begin{algorithmic}[1]
\State{$\calD$ contains $n$ points in $\bbR^d$}
\State{$f(\calD)$ contains values of $f(\bx_j)$ for all $\bx_j\in\calD$}
\State{$\bq \in CH(\calD)$ is an interpolation point}
\State{Set $\tilde\bx\leftarrow \texttt{argmin}_{\bx_j\in\calD}\|\|\bx_j-\bq\|\|$} \Comment{\textit{$\tilde\bx$ is the closest point in $\calD$ to $\bq$}}
\State{Find a $d$-simplex ${\calS}$ in $DT(\calD)$ incident to $\tilde\bx$} 

\While{$\bq \not\in {\calS}$}
\State{Select the facet ${\cal F}$ of ${\calS}$ from which $\bq$ is visible}
\State{Complete a new $d$-simplex ${\calS}^*$ from the facet ${\cal F}$}
\State{Update ${\calS} \leftarrow {\calS}^*$}
\EndWhile
\State{Since the loop has terminated, $\bq\in {\calS}$}\\
\Return{${\hat f}_{DT}(\bq)$} \Comment{\textit{Computed according to (\ref{eq:fDTx}), using $f(\calD)$}}
\end{algorithmic}
\end{algorithm}

Note that the cost to build the initial seed simplex is ${\cal O}(nd^3)$, where $n = \|\calD\|$, and the cost to compute each subsequent simplex is ${\cal O}(nd^2)$.
Therefore, the total cost of~\cref{alg:DS} is ${\cal O}(nd^3 + nd^2p)$, where $p$ is the number of flips required.
Empirically, for uniformly distributed $\calD$, $p$ tends to be a super-linear but sub-quadratic function of $n$ and independent of $d$~\cite{chang2018polynomial}.
In particular, typically $p\gg d$ so that Algorithm 1 is effectively ${\cal O}(nd^2p)$.

While~\cref{alg:DS} only allows interpolation points in the geometric sense, 
it can also extrapolate a value $f(\by)$ for a point $\by$ outside $CH(\calD)$.
For this, the point $\by$ is projected onto the convex hull by solving a quadratic programming problem, and the resulting projection ${\hat \by}\in CH(\calD)$ is interpolated on the face of a simplex in $DT(\calD)$ via~\cref{alg:DS}.
In this work, we detect and flag when a query point corresponds to geometric extrapolation for a given sample set, but intentionally set up experiments so that extrapolation does not occur.





\section{Computing the Delaunay density diagnostic} 
\label{sec:Delaunay}

We now describe in detail how the Delaunay density diagnostic is computed.
As outlined in~\cref{sec:intro}, the user specifies the following:
a set of \textbf{query points} $\{\bq_i\}\subset\bbR^d$; a $d$-dimensional \textbf{bounding box} $\calB\subset\bbR^d$, such that  $\{\bq_i\}\subset \calB$; an initial set of \textbf{sample points} $\{\bx_j\}\subset\bbR^d$, drawn randomly from $\calB$; and an \textbf{upsampling growth factor} $b\in (1,2]$.
In addition, a stopping criterion should be specified, in the form of a maximum number of upsampling iterations and/or a maximum size for the set of sample points.
The algorithm is as follows:

\begin{algorithm}
\caption{Delaunay density diagnostic (for \texttt{MSD} rate)}
\label{alg:DD_msd}
\begin{algorithmic}[1]
\State{$k$=0}
\While{stopping criteria not met}
\State{$n_k\leftarrow\|\{\bx_j\}\|$}
\State{${\hat f}_k(\bq_i)\leftarrow$ results of~\cref{alg:DS} with $\calD=\{\bx_j\}$, for each $\bq\in\{\bq_i\}$}
\If{$k>0$}
\State{$\text{diff}_{k,i} \leftarrow ({\hat f}_k(\bq_i)-{\hat f}_{k-1}(\bq_i))$, for each $\bq\in\{\bq_i\}$} \label{algDD:line:diff}
\EndIf
\If{$k>1$}
\State{$r_k \leftarrow \log_{\,b} \left(\frac{\sqrt{\msdi{\text{diff}_{k-1,i}}}}{\sqrt{\msdi{\text{diff}_{k,i}}}}\right)$}
\label{algDD:line:rate}
\EndIf
\State{$\ell\leftarrow \texttt{round}\left(\left[ b~n_k^{1/d}- (b-1)\right]^d -n_k, ~~0\right)$} \label{algDD:line:upsample}
\State{Generate $\ell$ points from $\calB$ randomly and add them to the collection $\{\bx_j\}$}
\State{$k\leftarrow k+1$}
\EndWhile\\
\Return{$\{n_k,~r_k\}$}
\end{algorithmic}
\end{algorithm}

Each output $(n_k,r_k)$ from~\cref{alg:DD_msd} is an estimate of the rate ($r_k$) at which piecewise linear interpolants are converging to $f$ at the query points $\{\bq_i\}$ when $n_k$ samples are drawn from $\calB$.
Since~\cref{alg:DD_msd} involves random sampling, we can run it multiple times with different initial random seeds to generate a distribution of outputs.
A collection of outputs $\{n_k,r_k\}_1$, $\{n_k,r_k\}_2$, $\ldots$ can then be used to assess the accuracy of the estimates and test sensitivity to the randomized aspect of the algorithm.

\paragraph{Computing rates in other norms}
To compute the rate in a different norm than mean squared difference (\texttt{MSD}), only a few lines from \cref{alg:DD_msd} need to be modified.  
An informative alternate rate to consider is the mean squared difference of the \textit{gradients} of the interpolants, presented in \cref{alg:DD_grad-msd}. 
We call this rate \texttt{grad-MSD}, for clarity, but note that it is essentially a discrete version of the Sobolev semi-norm $H^1$, just as the square root of $\texttt{MSD}$ is a discrete version of the $L^2$ norm.

\begin{algorithm}
\caption{Delaunay density diagnostic (for \texttt{grad-MSD} rate)}
\label{alg:DD_grad-msd}
\begin{algorithmic}[1]
\State{$\cdots$}
\setalglineno{4}
\State{$\nabla{\hat f}_k(\bq_i)\leftarrow$ results of~\cref{alg:DS} with $\calD=\{\bx_j\}$, for each $\bq\in\{\bq_i\}$} \label{alggdd:line:grad}
\If{$k>0$}
\State{$\text{diff}_{k,i} \leftarrow (\nabla{\hat f}_k(\bq_i)-\nabla{\hat f}_{k-1}(\bq_i))$, for each $\bq\in\{\bq_i\}$ } 
\EndIf
\If{$k>1$}
\State{$r_k \leftarrow \log_{\,b} \left(\frac{\texttt{norm}({\text{diff}_{k-1,i}})}{\texttt{norm}({\text{diff}_{k,i}})}\right)$}  \label{alggDD:line:norm}
\EndIf
\setalglineno{14}
\State{$\cdots$} \\
\Return{$\{n_k,~r_k\}$}
\end{algorithmic}
\end{algorithm}
We briefly comment on two key changes from \cref{alg:DD_msd} to \cref{alg:DD_grad-msd}. First, the calculation of $\nabla{\hat f}_k(\bq_i)$ in line~\ref{alggdd:line:grad} can be done accurately and robustly since $\hat f_k$ is piecewise linear.
Details of this computation are provided in~\cref{app:grad}.
Second, we interpret the term \texttt{norm} used in line~\ref{alggDD:line:norm} as the regular Euclidean norm for vectors in $\bbR^d$, which is available as \texttt{numpy.linalg.norm} in Python.

We have equivalent claims to \cref{claim:rate2} and \cref{claim:rate0} for the \texttt{grad-MSD} rate; the only change is the expected value for $r_k$ in each case, namely, 1 replaces 2 and -1 replaces 0.

\begin{claim}
\label{claim:gradrate1}
If the set $\{(\bx_j,f(\bx_j)\}$ used to define $\hat f_k$ contains sufficient information to reconstruct the geometric features present in $f$ over the set of query points $\{\bq_i\}$ then the \texttt{grad-MSD} rate $r_k$ will be approximately $1$.
\end{claim}

\begin{claim}
\label{claim:gradrateminus1}
If the set $\{(\bx_j,f(\bx_j)\}$ used to define $\hat f_k$ cannot distinguish geometric features from random noise over the set of query points $\{\bq_i\}$ then the \texttt{grad-MSD} rate $r_k$ will be approximately $-1$.
\end{claim}

Other norms could be considered as well.  
For instance, computing rates in~$L^p$ for $p\in[1,\infty)$ could be assessed by replacing $\texttt{MSD}$ with ``mean $p$th power differences'' and the square roots with $p$th roots in line \ref{algDD:line:rate} from \cref{alg:DD_msd}.
The $L^\infty$ rate could be similarly accomodated.
Likewise, convergence in the Sobolev norm $W^{1,p}$ can be attained by using the $L^p$ norm instead of the $L^2$ norm for the term \texttt{norm} used in line~\ref{alggDD:line:norm} of \cref{alg:DD_grad-msd}.
We leave study of these possibilities for future work.

It is important to note that the point of using a linear interpolant is not because it is necessarily the best fit for the data. Indeed, there are many algorithms that will likely fit the data ``better,'' in the sense of having lower generalization error.  We use a linear interpolant here because it has a known convergence rate for smooth data, which is required to serve as a diagnostic. Most algorithms, especially ones requiring any type of iterative fitting or training procedure, do not have known convergence rates and thus cannot be used as intended here.  Further, since we are only interested in the convergence \emph{rate} from successive approximations, the quality of the interpolant as a surrogate function is irrelevant. We choose the Delaunay linear interpolant because of its unique definition for arbitrary datasets and for the speed with which \texttt{DelaunaySparse} is able to work in moderate dimensions.

\section{Numerical results} 
\label{sec:numerics}

We have implemented and tested the feasibility of the Delaunay density diagnostic in various scenarios for data in $\bbR^d$ from $d=2$ to $d=5$.  
After collecting the data, we convert the number of samples at the $k$th step, $n_k$, to average sample spacing by
\begin{equation}
\label{eq:as-def}
\as := 
\text{average sample spacing} := \frac{L}{n_k^{1/d}},
\end{equation}
where $L$ is the average side length of the bounding box $\calB$.
For simplicity, we have always taken $\calB$ to be a $d$-dimensional cube in $\bbR^d$; the definition of $L$ could be modified if a more complicated choice of $\calB$ was desired.
The values of average sample spacing are used for the horizontal axes in our figures.
Thus, the smallest average sample spacing in each plot corresponds to the \textit{largest} $n_k$ value attained before the stopping criteria.

To assess the sensitivity of~\cref{alg:DD_msd} and~\cref{alg:DD_grad-msd} to the location of points $\{\bx_j\}$, we use a seed to initialize a random number generator.
At the start of the code, we set \texttt{rng=numpy.random.default\_rng(globalseed)},
where \texttt{globalseed} is a user-specified integer; we then use $\texttt{rng}$ to generate samples $\{\bx_j\}$ subsequently in the code.
After running the code with multiple seeds, we report the ``mean rate'' computed over all seeds, as well as the inter-quartile and inter-decile range for the computed rates, thereby providing an estimate of uncertainty in the computation.\\

\begin{figure}[t]
	\centering
	\footnotesize
	\begin{tabular}{cc}
		\includegraphics[width=.35\textwidth]{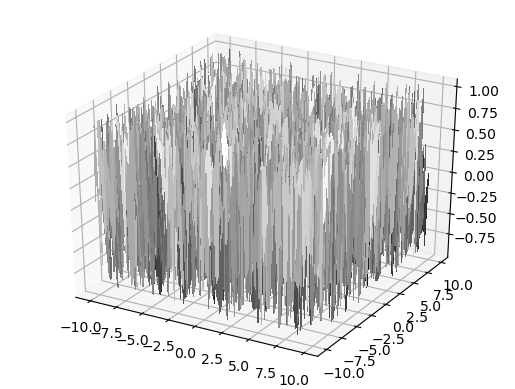} &
		\includegraphics[width=.35\textwidth]{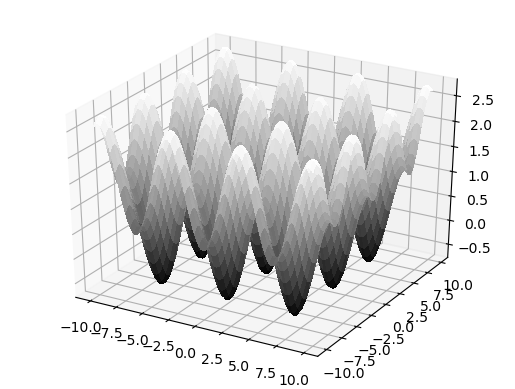} \\
		\includegraphics[width=.35\textwidth]{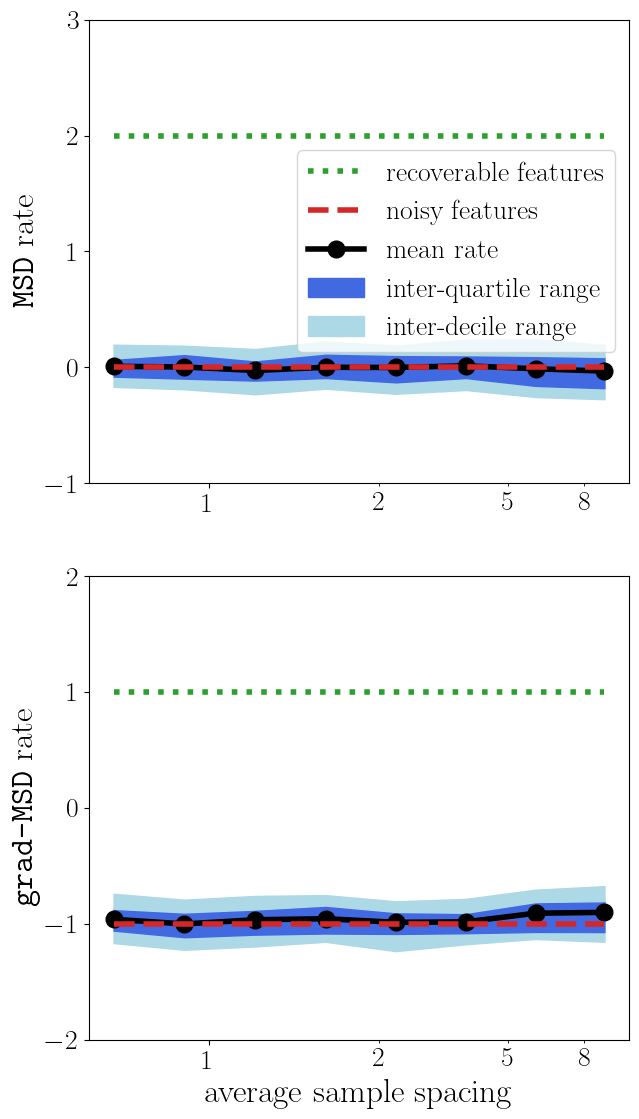} &
		\includegraphics[width=.35\textwidth]{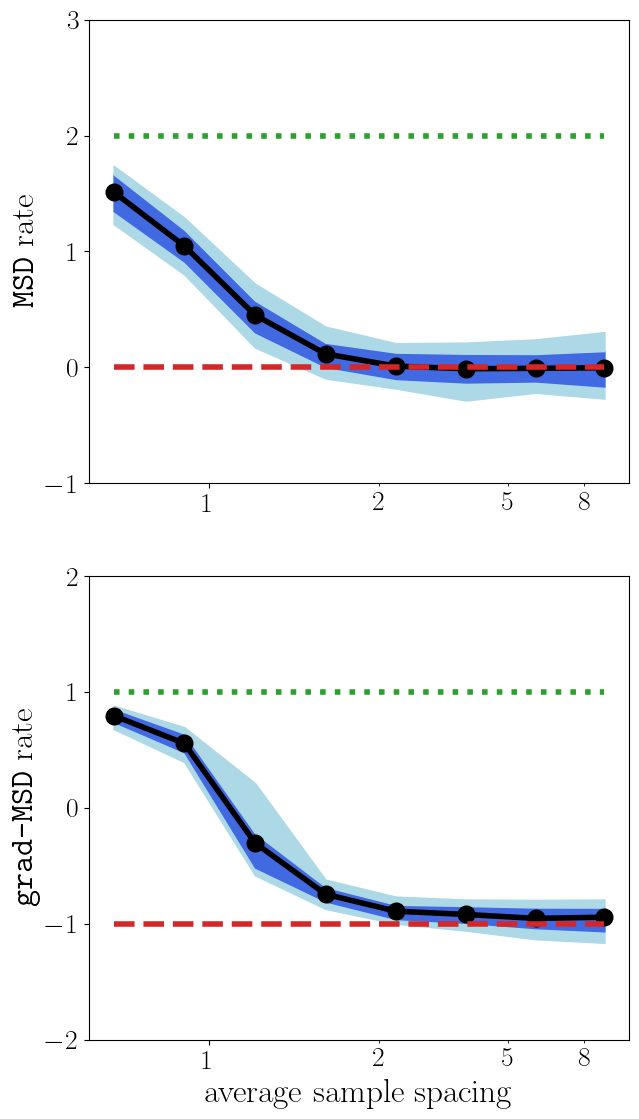} \\
		        (a) & (b) \\
	\end{tabular}
	\caption{We validate~\cref{alg:DD_msd} and~\cref{alg:DD_grad-msd} by assessing the computed rate over a range of average sample spacings, for $f$ shown in the top row.
	The mean rate (black dot series) shows the average of the computed rate over 100  trials with different random initial seeds. The inter-quartile (between $25^\text{th}$ and $75^\text{th}$ percentiles) and inter-decile (between $10^\text{th}$ and $90^\text{th}$ percentiles) ranges are shown in dark blue and light blue bands, respectively.
	(a) Pure noise is consistently detected as having the ``noisy features'' rate (0 for \texttt{MSD}, $-1$ for \texttt{grad-MSD}).
	(b) Fine scale features are only recoverable if average sample spacing is small enough.
}
	\label{fig:validation1}
	\vspace{-.2in}
\end{figure}

\begin{figure}[t]
	\begin{tabular}{cc}		        
		\includegraphics[width=.35\textwidth]{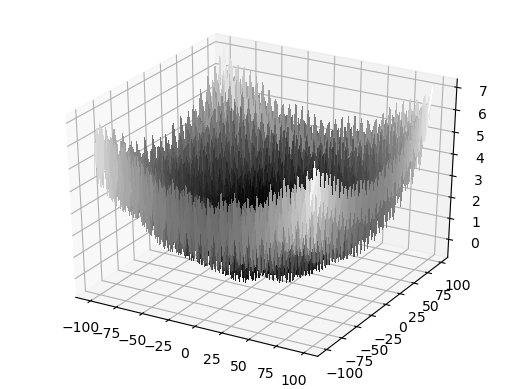} &
		\includegraphics[width=.35\textwidth]{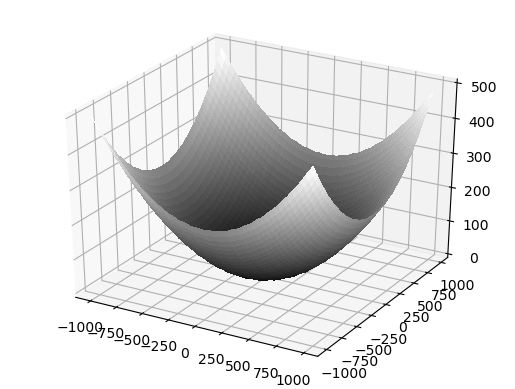} \\
		\includegraphics[width=.35\textwidth]{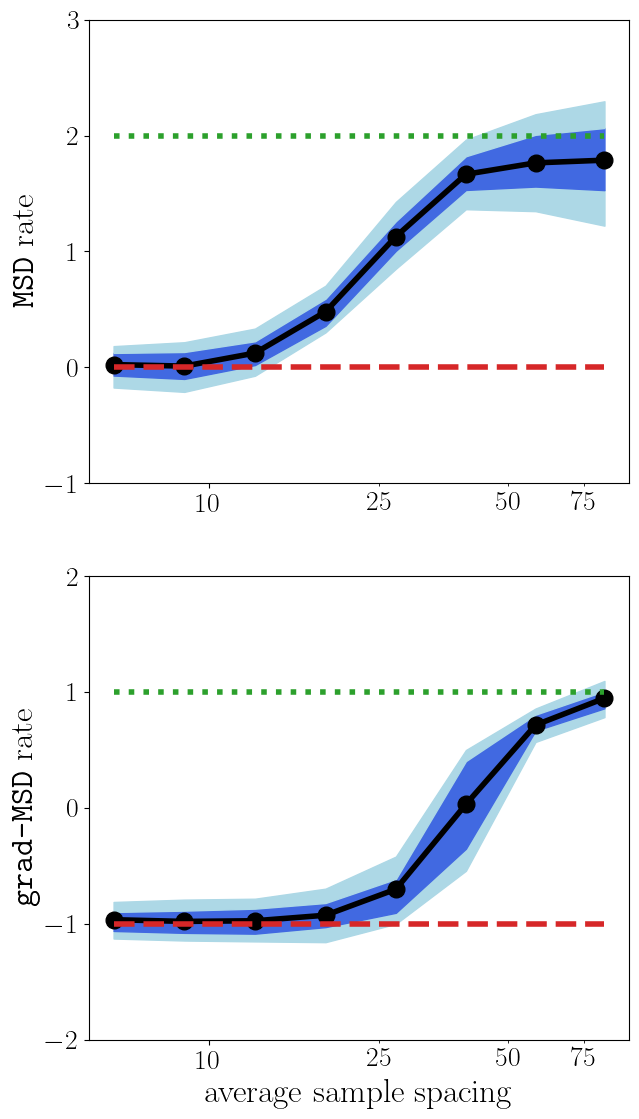} &
		\includegraphics[width=.35\textwidth]{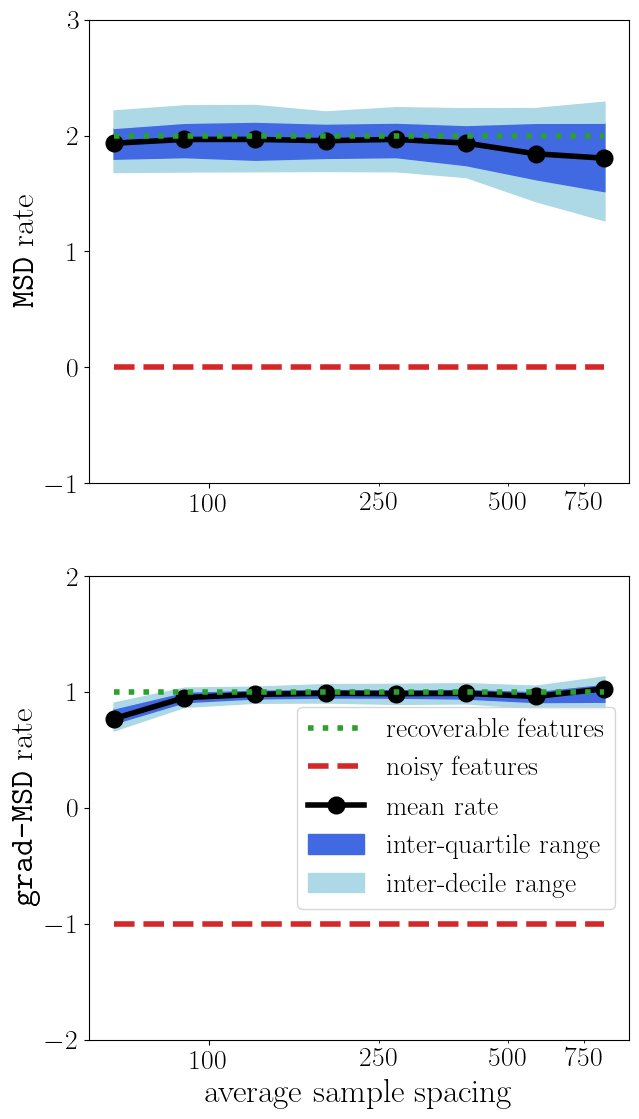} \\
 (a) & (b) \\
	\end{tabular}
	\caption{
	Additional experiments as in~\cref{fig:validation1}.
	(a) Insufficiently sampled fine scale features are detected as noise; the large scale quadratic feature is recoverable with larger average sample spacing.  
	(b) Smooth variation (at the scale of inquiry) is consistently detected as having the ``recoverable features'' rate (2 for \texttt{MSD}, 1 for \texttt{grad-MSD}).	}
	\label{fig:validation2}
	\vspace{-.2in}
\end{figure}

\paragraph{Distinguishing features from noise}
We begin by demonstrating that the rate computations can reliably distinguish features from noise for analytic functions $f:\bbR^2\raw\bbR$, shown in~\cref{fig:validation1,fig:validation2}.
For these examples, we set the query points $\{\bq_i\}$ to be a $20\times 20$ uniformly spaced grid including the corners of $[-10,10]^2$ for (a) and (b), $[-100,100]^2$ for (c), and $[-1000,1000]^2$ for (d).
The bounding box $\calB$ is taken to be a square centered at the origin with side length 25\% longer than the query point lattice, e.g.,~$\calB=[-12.5,12.5]^2$ for (a).
The upsampling growth factor $b$ is 1.4641.
We initialize~\cref{alg:DD_msd} and~\cref{alg:DD_grad-msd} with $\|\{\bx_j\}\|=9$ and use a stopping criterion of $\|\{\bx_j\}\|>200,000$.

When calling DelaunaySparse (\cref{alg:DS}), we compute and pass values of $f(\{\bx_j\})$ for the current collection of sample points $\{\bx_j\}$.
In (a), we assign $f(\bx_j)$ to be a random number drawn uniformly from $[-1,1]$.
In (b)$-$(d), we assign $f(\bx_j):=g_2(\bx_j)$, where $g_2$ is the Griewank function~\cite{G1981} on $\bbR^2$, given by:
\begin{equation}
\label{eq:grwk}
g_d(x_1,\ldots,x_d) := \sum_{i=1}^d\frac{x_i^2}{4000}-\prod_{i=1}^d\cos\left(\frac{x_i}{\sqrt i}\right) + 1.
\end{equation}

The top rows of~\cref{fig:validation1,fig:validation2} show a visualization of $f$ for each case, over a domain matching the region defined by the query points.
The middle and bottom rows show the \texttt{MSD} and \texttt{grad-MSD} rates, respectively, as a function of average sample spacing.
Here, we have filtered out $(n_k,r_k)$ outputs from the code with $n_k<500$, i.e.~larger average sample spacings.
These outputs had larger variations in $r_k$ values due to the small number of points involved, which distracted from the success of the method for larger $n_k$ values.  
We will discuss this issue further later.
Thus, in~\cref{fig:validation1,fig:validation2}, the average sample spacing values correspond to $n_k\in\{
173832,
~81275,
~38039,
~17830, 
~8376, 
~3947, 
~1869,   
~891
\}$, 
from left to right, in each graph.
The captions of~\cref{fig:validation1,fig:validation2} describe high level takeaway messages from each experiment.\\





\renewcommand{\arraystretch}{2}
\begin{figure}[t]
	\centering
	\footnotesize
	\begin{tabular}{cc}
		\includegraphics[width=.5\textwidth]{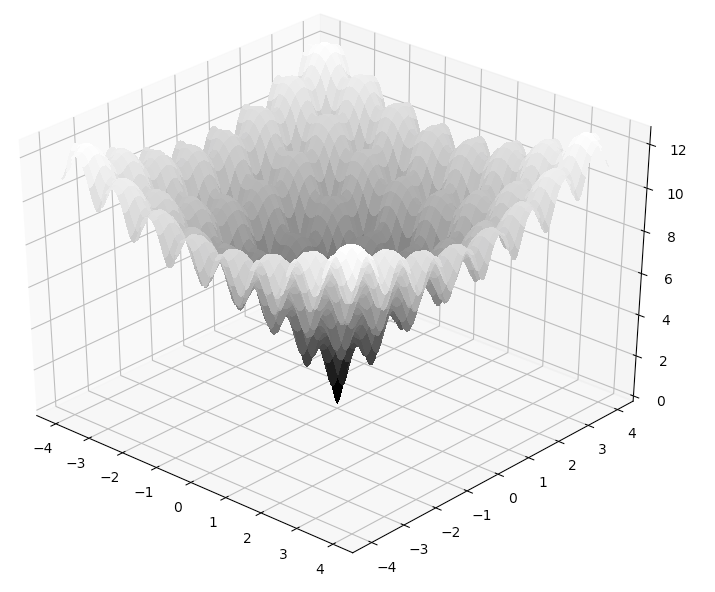} 
	&	\includegraphics[width=.25\textwidth]{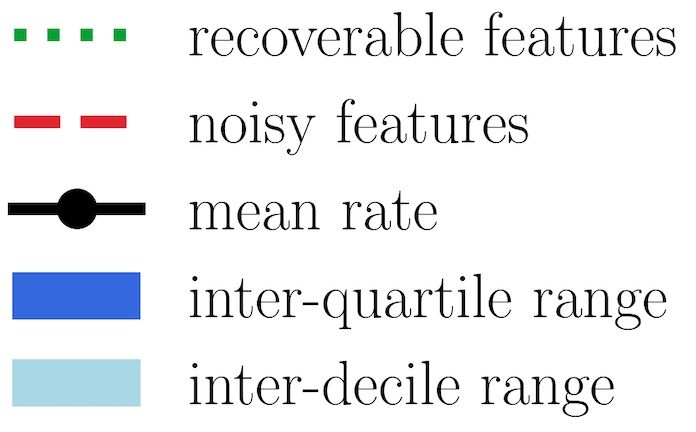}
	\end{tabular}

	\caption{\textit{Left}: A basic visualization of the Ackley function on $[-4,4]^2$.  \textit{Right}: Legend for the plots in~\cref{fig:ackley2}.}
	\label{fig:ackley1}

\end{figure}

\renewcommand{\arraystretch}{2}
\begin{figure}[t]
	\centering
	\footnotesize
	\begin{tabular}{ccc}
		\includegraphics[width=.30\textwidth]{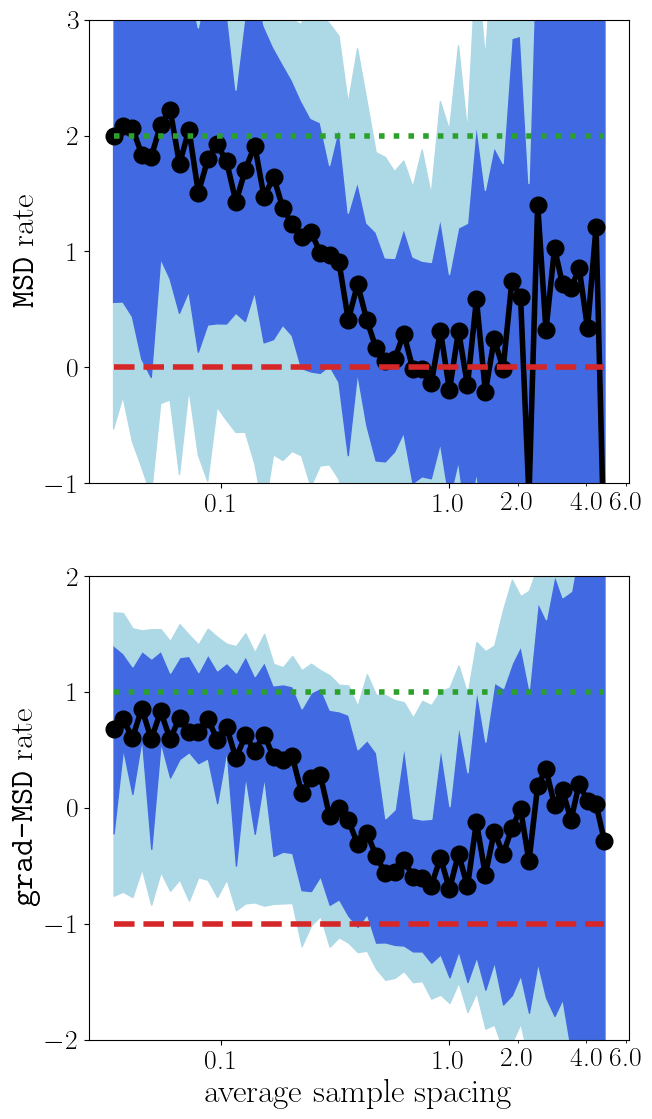} &
		\includegraphics[width=.30\textwidth]{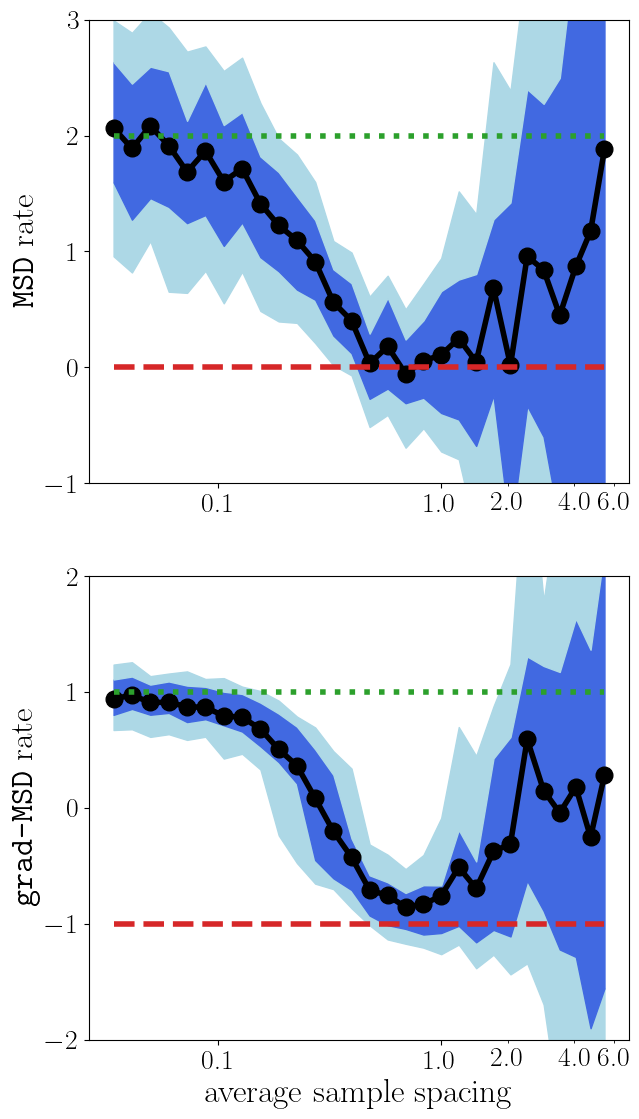} &
		\includegraphics[width=.30\textwidth]{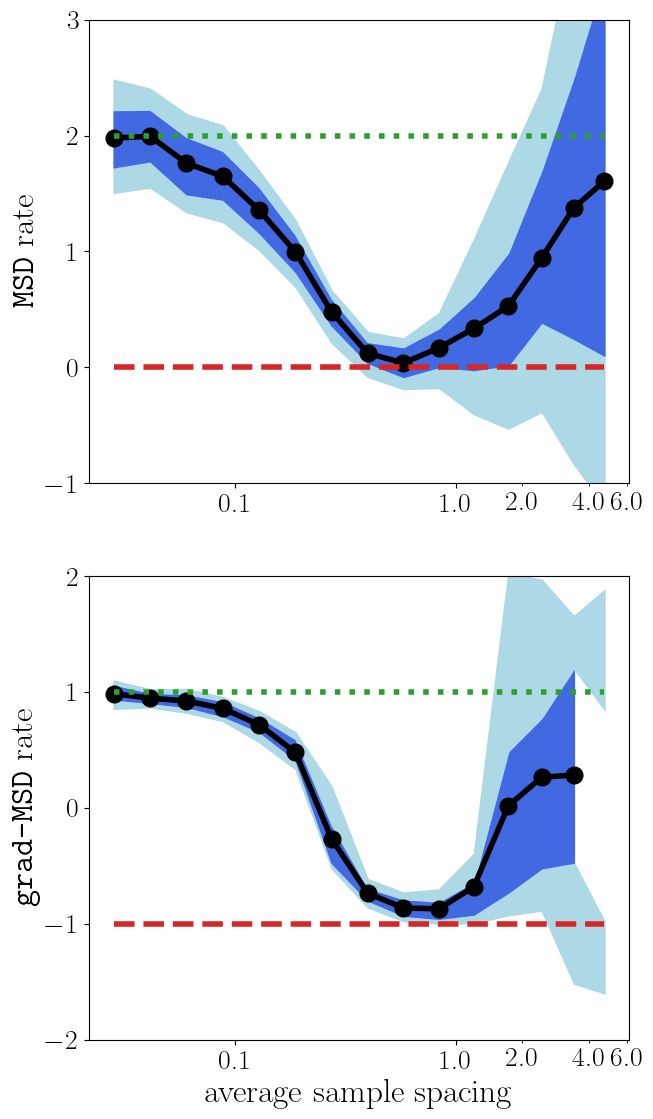} \\
		 $b=1.1$ & $b=(1.1)^2=1.21$ & $b=(1.1)^4=1.4641$ \\
		 	\vspace{-.2in}
	\end{tabular}

	\caption{For the Ackley function on $\bbR^2$---see~\cref{fig:ackley1}---we examine the effect of changing the upsampling rate $b$ on the rates computed by~\cref{alg:DD_msd} and~\cref{alg:DD_grad-msd}.  Using the values of $b$ indicated, we confirm that a larger $b$ value corresponds to a larger step size (in the horizontal axis) and a smaller variation in the computed rate, as evidenced by the narrower inter-quartile and inter-decile ranges as $b$ increases.  For larger average sample spacings, the small number of samples is the cause of the increased variation.  In each case, the mean rate has the same trend, reflecting the multi-scale nature of the Ackley function.}
	\label{fig:ackley2}
	\vspace{-.2in}
\end{figure}

\paragraph{Effect of the upsampling factor $b$}
We next examine the effect of varying $b$, the parameter that controls the rate at which samples are added during each iteration (line~\ref{algDD:line:upsample} from Algorithm~\ref{alg:DD_msd}). 
For this experiment, we fix $f$ to be the Ackley function~\cite{A2012}, visualized in~\cref{fig:ackley1}, given by the formula
\[f_{\text{(Ackley)}}(x,y) := -20~\exp\left(-0.2 \sqrt{\frac{x^2 + y^2}{2}}\right) - \exp\left(\frac{\cos(2\pi x) + \cos(2\pi y)}{2} \right) + 20 + e.\]
Like the Griewank function, the Ackley function is commonly used to test optimization algorithms due to its dense collection of local extrema at a fine scale.
We again fix the query points $\{\bq_i\}$ to be a $20\times 20$ uniformly spaced grid including the corners of $[-10,10]^2$ and set $\calB=[-12.5,12.5]^2$.   
We initialize each trial with $\|\{\bx_j\}\|=9$ and use a stopping criterion of $\|\{\bx_j\}\|>1,000,000$.

We carried out 100 trials for three different values of the upsampling growth factor $b$: 1.1, 1.21, and 1.4641, shown in~\cref{fig:ackley2}.
As in the previous example, we filter out results with very small $n_k$ values (in this case $n_k<20$) due to the wide variation in $r_k$ values, as evidenced by the larger error bars on the right side of each graph.
Regardless, the trend picked up by the mean rate in each experiment is quite clear: the Ackley function has many small oscillations with a period of about 1 unit in each coordinate, so we can fully resolve the function as long as the sample density has $\as\ll 1$.
\cref{fig:ackley2} indicates further that we can expect to fully recover features when $\as\lesssim 0.1$. Additionally, when looking at the Ackley function over larger scales (several units), the oscillations are less important to the value of the function than the global bowl-shaped trend. \cref{fig:ackley2} demonstrates this as well: for $\as>5$, we can fully recover the large-scale behavior of the function. In between these regimes, our diagnostic is detecting ``features'' (the oscillations) at about $\as\sim 0.8$, roughly matching the scale of the oscillations (the period of the cosine terms, i.e.~1 unit).\\

\begin{figure}[t]
	\centering
	\footnotesize
	\begin{tabular}{ccc}
	\includegraphics[width=.30\textwidth]{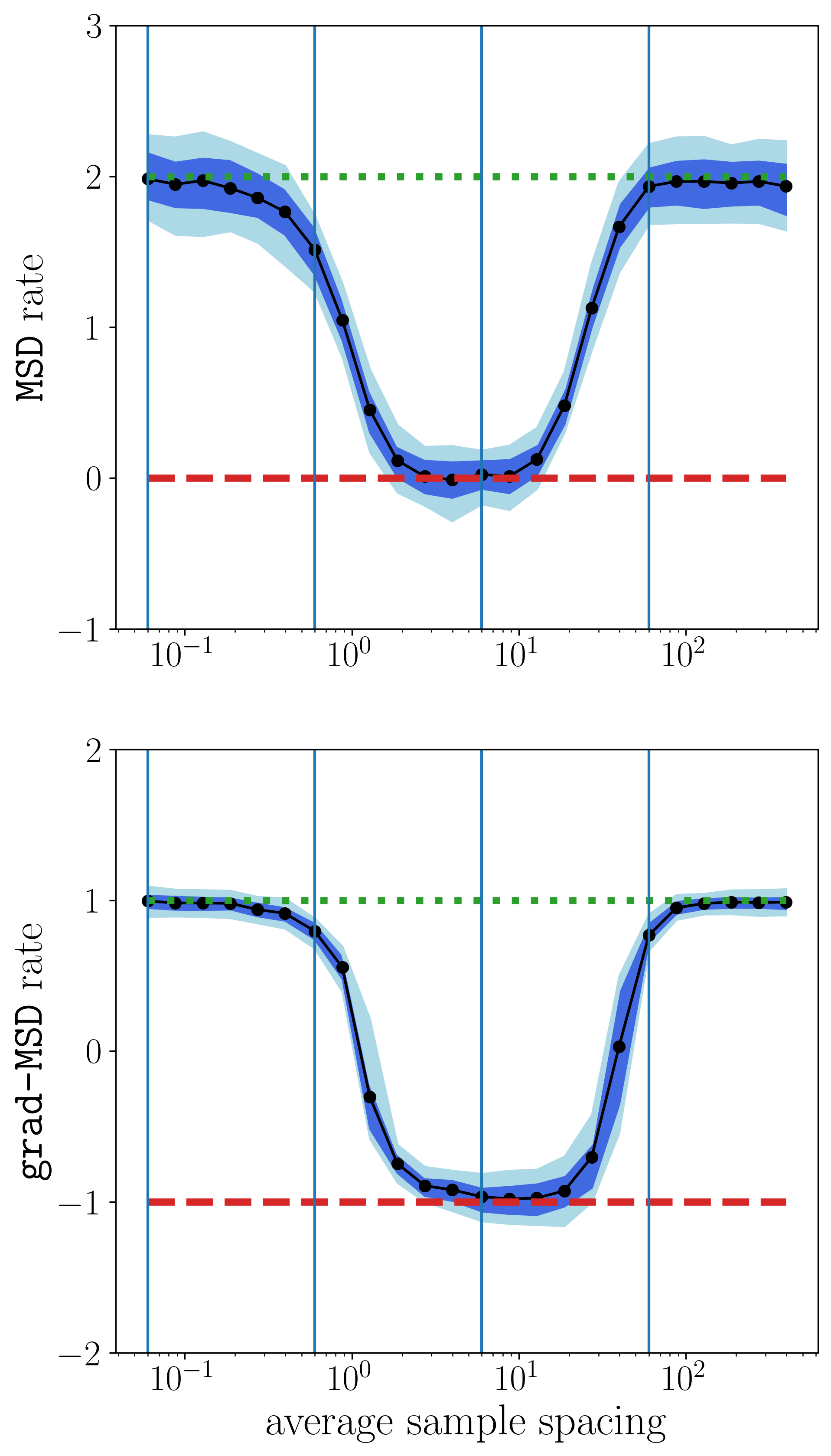} &
	\includegraphics[width=.30\textwidth]{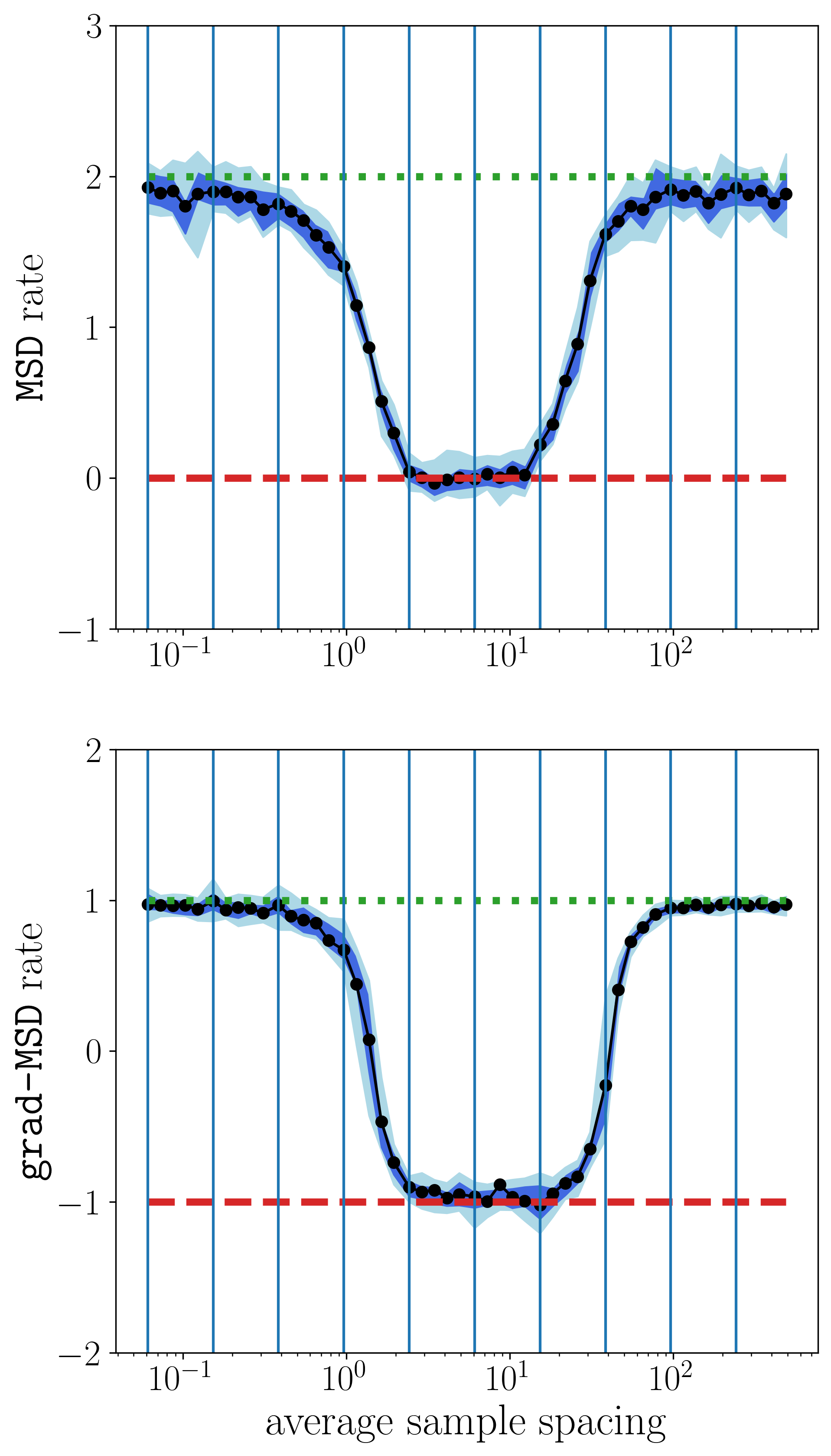} &
	\includegraphics[width=.30\textwidth]{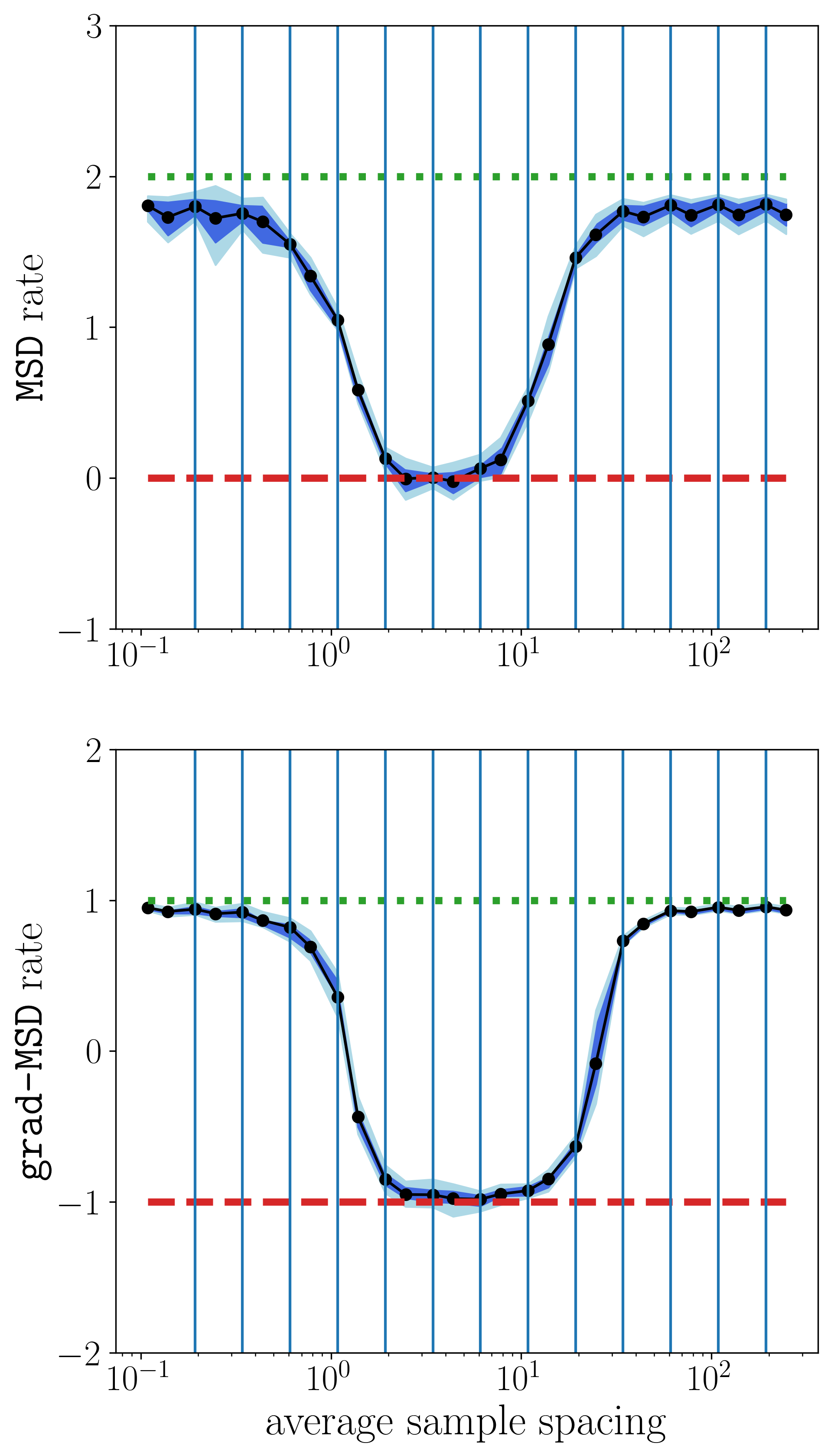} \\
	$f = g_2$ & $f = g_3$ & $f = g_4$ \\
	\end{tabular}
	\caption{We examine the effect of stepping up the input dimension for the Griewank function $g_d:\bbR^d\raw\bbR$, given in \cref{eq:grwk}.  We fix query points $\{\bq_i\}$ and a bounding box $\calB$, centered at the origin, then run distinct experiments by scaling both by pre-determined amounts.  Vertical blue bars separate the distinct scales of experiments, i.e., places where we adjusted the size of the bounding box. As the dimension increases, we use more distinct scales with fewer iterations per scale, demonstrating one approach for scaling with dimension. The computed rates have the same behavior in each dimension, further validating the method.}
	\label{fig:grwk2d3d4d}
	\vspace{-.2in}
\end{figure}

\paragraph{Experiments in higher dimensions}
While many Delaunay methods for computational geometry apply exclusively to data in $\bbR^2$ or $\bbR^3$, the Delaunay density diagnostic algorithms have no formal restriction on input dimension $d$.
We explore practical considerations of assessing data sets with dimension 4 or higher by a series of experiments with the Griewank functions, defined in (\ref{eq:grwk}).
While multiple parameters must be selected in order to run the code, each parameter has a clear geometric interpretation of how it affects the computational cost and accuracy of the result, as we will explain.
The results of the experiments are shown in~\cref{fig:grwk2d3d4d}; the parameters used are given in~\cref{tab:gwkdimparams} and explained below.

For $d=2$, 3, 4, we fix a uniformly spaced $d$-dimensional lattice of $Q^d$ query points, centered at the origin in $\bbR^d$.  
For $d=2$, $3$ we use $Q=20$ and for $d=4$ we use $Q=10$.
Note that $Q$ can be any positive integer ($Q=1$ corresponds to a single query point), however, the size of $Q^d$ will become a main driver of computational cost as $d$ increases.
Letting $M$ denote the side length of the query point lattice, and $L$ the side length of the bounding box to be used (as in~\cref{eq:as-def}), we define 
\[\qpdf := \text{query points dimension fraction} := \frac ML. \]
For $d=2$, $3$, we use $\qpdf=0.8$ and for $d=4$ we use $\qpdf=0.6$, with the bounding box always centered at the origin.
The purpose in decreasing $\qpdf$ for larger $d$ is to reduce the probability that a query point lies outside the convex hull of the samples $\{\bx_j\}$, a case we exclude from consideration in this work, per the discussion at the end of \cref{subsec:DS}.
    
\renewcommand{\arraystretch}{2}

We select $b$ as indicated in~\cref{tab:gwkdimparams}. 
The effect of $b$ on accuracy was discussed above and in~\cref{fig:ackley2}.
We initialize each trial with $\|\{\bx_j\}\|=3000$ for $d=3$ and 5000 for $d=4$ and set a stopping criteria of $\|\{\bx_j\}\|=100,000$ or $200,000$.
After collecting results, we filter out small $n_k$ values to produce the min and max values indicated in~\cref{tab:gwkdimparams}.
From these values, we compute $\as$ and then determine the number and extent of scaling needed so that $\as$ will have range $\approx [10^{-1},10^{2.5}]$ across all experiments.  This last point is key. There are not enough computational resources to support exploring the Griewank function at scales ranging from $10^{-1}$ to $10^{2.5}$ using a single bounding box (too many points are required to reach the smallest scales). So instead, we scale down the bounding box when we reach our limit of points (100k-200k) and then glue these different scales together by matching $\as$.

\begin{table}[h]
\centering
\begin{tabular}{crrrrrcc}
dim & $b$ & $\|\{\bq_i\}\|$  & \qpdf &  $\min n_k$ &  $\max n_k$ & \# scales  & \# trials\\
\hline
2 & 1.4641 & 400    & 0.8 & 3,947  & 173,832 & 4  & 100 \\
3 & 1.2    & 8,000  & 0.8 & 8,400  & 69,321  & 10 & 25  \\
4 & 1.3    & 10,000 & 0.6 & 33,423 & 89,109  & 14 & 25
\end{tabular}
\caption{Parameters used to generate results shown in~\cref{fig:grwk2d3d4d}.}
\label{tab:gwkdimparams}
\end{table}
 
We use high performance computing resources at Lawrence Livermore National Laboratory to compute all of our examples. 
Each trial was run on a single node with 64 GB memory, consisting of 16 cores.
We call the built-in ``level 1 parallelism'' feature of \texttt{DelaunaySparse}, which exploits a speedup strategy on the loop over query points; this approach temporarily stores discovered parts of the Delaunay mesh structure to accelerate subsequent interpolation queries.
With these resources and parameter choices, the wall clock time per trial was 45 seconds, 8.5 minutes, and 11.1 minutes in $d=2$, 3, and 4, respectively.
Thus, the cumulative compute time to produce~\cref{fig:grwk2d3d4d} was approximately 5 hours, 35 hours, and 65 hours for $d=2$, 3, and 4, respectively. 

The reported run times indicate the feasibility of the algorithm in higher dimensions, not an actual assessment of its scalability with respect to dimension.
We did not attempt to optimize the code for speed.
In practice, the dimension of the data is fixed by the application context.
A user would only need to adjust the other parameters---i.e.,~the column labels of~\cref{tab:gwkdimparams}---to ensure the error bounds met their requirements and the compute time fit within their resource capabilities.

\paragraph{Fixed datasets in $\bbR^d$ for small values of $d$}
We now demonstrate how the algorithm can be modified to assess sampling density of static, existing datasets.
At a high level, the only major change required to employ~\cref{alg:DD_msd} or \cref{alg:DD_grad-msd} on a static dataset is a change to the random point generation process.
No bounding box is used.
Instead, the index for the static dataset is randomly shuffled and the initial collection of sample points $\{\bx_j\}$ is defined to be the first $n_0$ points indicated by the shuffled index.
Subsequent additions to $\{\bx_j\}$ are attained by including the next $\ell$ points according to the shuffled index. 
This process emulates the random selection of points, avoids drawing duplicates from the dataset, and is limited by the size of the dataset.
We summarize this modification in~\cref{alg:DD_static}.

\begin{algorithm}
\caption{Delaunay density diagnostic (for static data sets)}
\label{alg:DD_static}
\begin{algorithmic}[1]
\State{$\cdots$}
\setalglineno{12}
\State{Select the next $\ell$ points from the static data with shuffled index. 
If $\ell$ points are not available, \textbf{break}.}
\setalglineno{14}
\State{$\cdots$} \\
\Return{$\{n_k,~r_k\}$}
\end{algorithmic}
\end{algorithm}

\begin{figure}[t]
	\centering
	\footnotesize
	\begin{tabular}{ccc}
		  \multirow{1}{*}{\includegraphics[width=.30\textwidth]{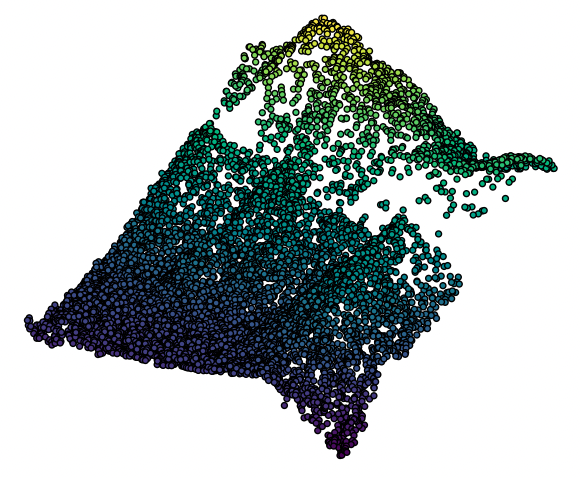}}
		& \multirow{2}{*}{\includegraphics[width=.30\textwidth]{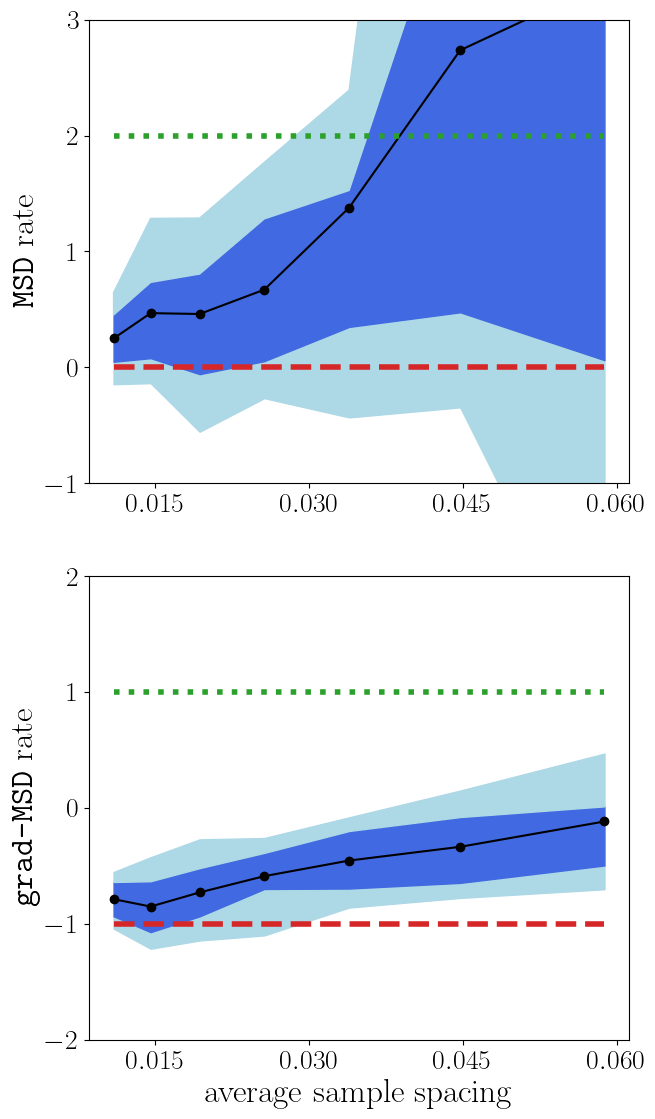}} 
		& \multirow{1}{*}{\includegraphics[width=.30\textwidth]{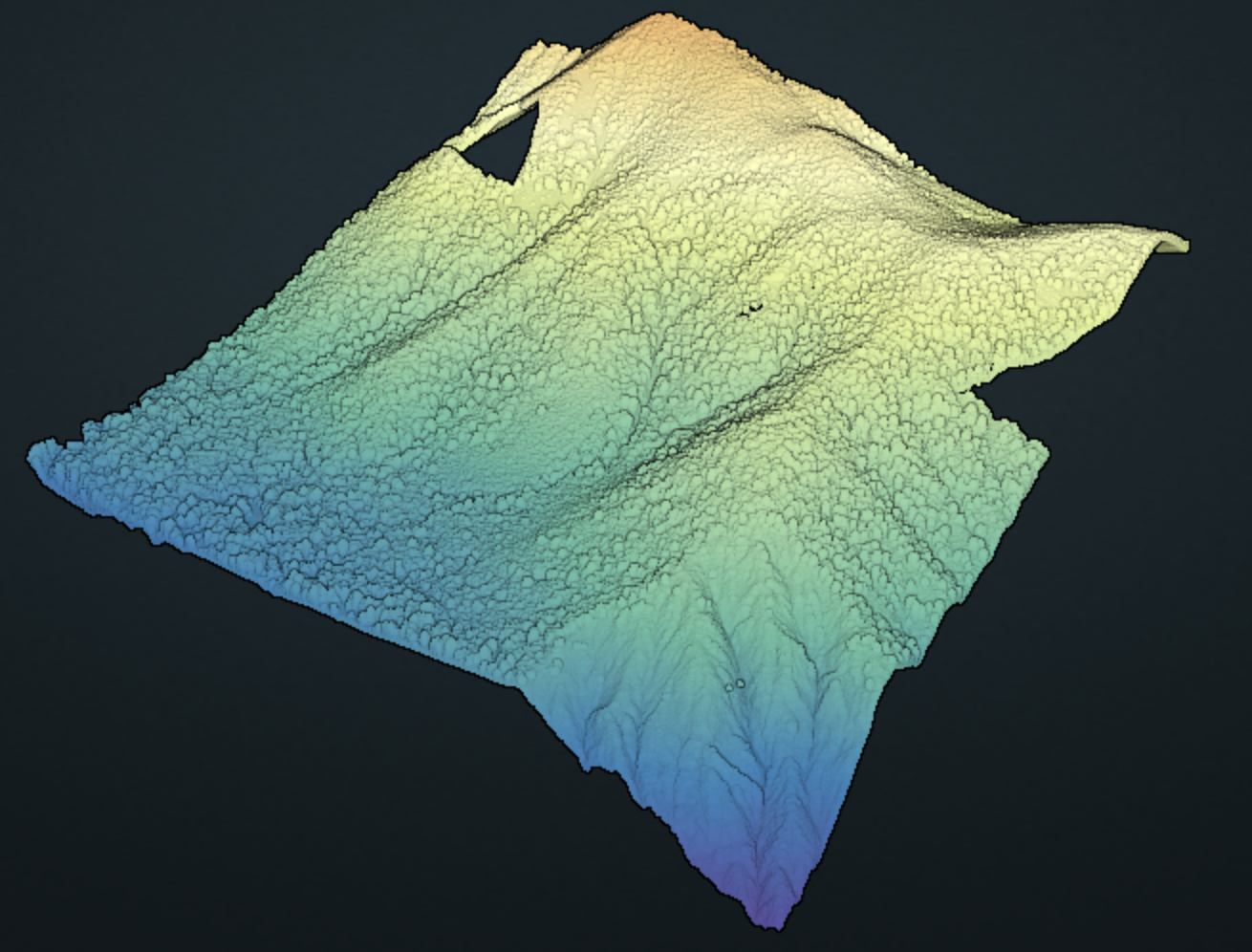}}\\
		& \\
		& \\
		& \\
		& \\
		& \\
		& \\
		\includegraphics[width=.30\textwidth]{figs/density_profile_legend.png} 
		& & \includegraphics[width=.30\textwidth]{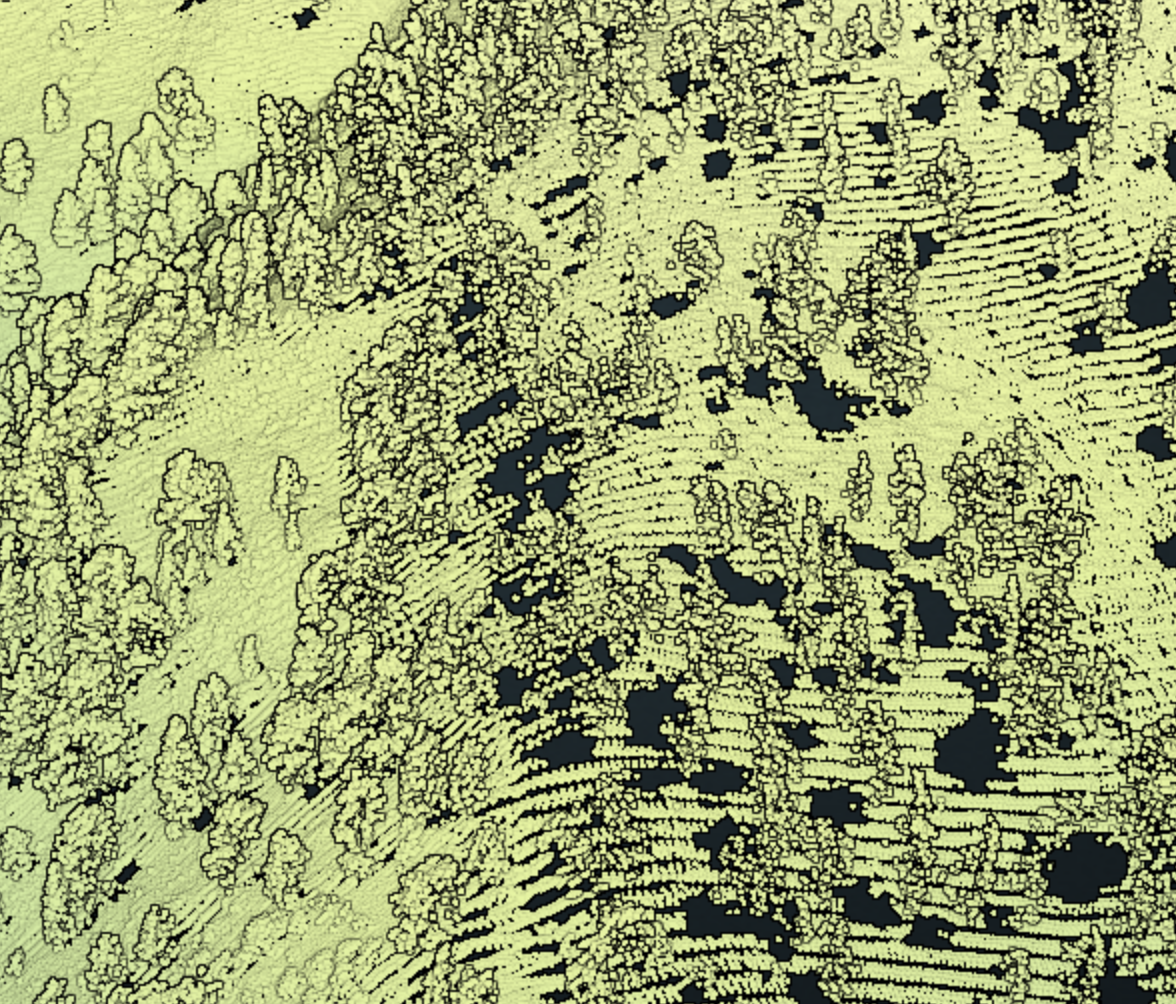} \\
	\end{tabular}
	\caption{We test~\cref{alg:DD_static} on a static dataset of 10,000 points, collected from a lidar sampling of a mountain range, retrieved from OpenTopography~\cite{C2016topo} (left, top).  The computed rates (middle) indicate that small scale features are present in the data set but not fully resolved.
These findings are consistent with visualizations of the original twenty million point cloud data from which the sample was drawn, indicating both coarse features at the scale of mountains (right, top) and fine features at the scale of trees (right, bottom). }
	\label{fig:topography}
\end{figure}

We have validated~\cref{alg:DD_static} on a variety of datasets, including datasets from the UCI Machine Learning Repository~\cite{UCIMLrepo} as well as fixed datasets of different sizes and dimensions generated from analytical functions.
We describe one such validation study in detail, using a dataset from OpenTopography~\cite{C2016topo}.
The full dataset is twenty million $x$, $y$, $z$ coordinates collected via lidar scanning over a roughly 4 $\text{km}^2$ region of a the San Bernadino Mountains in California.
A simple point cloud rendering reveals that there are large scale ``mountain size'' features as well as small scale ``tree size'' features present in the dataset; see Figure~\ref{fig:topography}, right column. 
We randomly sample ten thousand points from the full dataset and save this as a fixed dataset for testing~\cref{alg:DD_static}.

The topography dataset is treated as samples from a function $f:(x,y)\mapsto z$.
By applying~\cref{alg:DD_static} to this dataset, we are seeking to determine the scale of features present in the ten thousand point sample of $f$.
From other experiments, we find that setting $\max n_k$ to the size of the dataset is acceptable due to the random nature of subsampling and that setting $\min n_k$ to $1\%$ of the dataset size often yields informative results.  
Accordingly, we set $\min n_k = 100$ and $\max n_k = 10,000$.
We set $b=1.333521432163324$, derived using a heuristic to produce a handful of computed rates (details on the heuristic appear in the code~\cite{dddrepo}).
We construct a grid of 1600 query points over the 10th--90th percentiles of the $x$ and $y$ coordinate data, which corresponds to $\qpdf=0.8$, since the samples are close to uniformly distributed over a rectangular region.
We set \# seeds = 100.
Under these settings, each seed took a few seconds to compute on a MacBook Pro.

In Figure~\ref{fig:topography}, middle, we show the result of the experiment just described.
The rates display trends similar to the Griewank function on $[-100,100]^2$, as studied in~\cref{fig:validation2}(a).
Notice that both rates trend upward as the sample spacing increases, from ``noisy'' toward ``recoverable''.
Accordingly, we may infer that the sample encompasses some non-trivial large-scale features (i.e.~the major contours of the mountains) but also has under-resolved small-scale features (i.e.~the trees, which register as noise in this dataset).
The larger sample spacing results do not coalesce around the recoverable features lines, since those rates are computed starting from $\min n_k$, which is only one hundred points.
A takeaway message from these findings is as follows: the ten thousand point sample captures large-scale features to moderate accuracy, but a regression or smoothing technique would be advisable to mitigate the effect of noise from under-resolved small-scale features.

\renewcommand{\arraystretch}{2}

In our final example, we apply~\cref{alg:DD_static} to two higher-dimensional static datasets in $\bbR^5$.
The datasets, denoted \texttt{BURN\_OFF} and \texttt{BURN\_ON}, consist of input--output pairs from HYDRA, a multi-physics simulation code developed at LLNL over the past twenty years that informs experiments at the National Ignition Facility~\cite{Metal2001}.
Both datasets consist of pairs $\{(\bx_j,f(\bx_j)\}$ where $\bx_j\in\bbR^8$ and $f(\bx_j)\in\bbR^1$ is a quantity of interest. 
Physicists familiar with the problem context have indicated that the response of $f$ is primarily dependent on only five of the eight inputs.
Thus, we filter out the three less-important input variables, reducing inputs to $\bx_j\in\bbR^5$.
The filtering process creates ``near-duplicates'' in the dataset, i.e., points $(\bx_m, f(\bx_m))$ and $(\bx_n, f(\bx_n)$ where $\bx_m\not=\bx_n$ but $\vn{\bx_m-\bx_n}_{L^2(\bbR^5)}<\delta$, with $\delta$ small enough that~\texttt{DelaunaySparse} marks them as identical.
We fix a tolerance level $\delta$, then identify clusters of points within the dataset such that each point in a cluster is within $\delta$ of some other point in the cluster.
For clusters with more than one point, we keep the mean of the points and values as a ``new'' data point and discard the points defining the cluster.

By the filtering process, the \texttt{BURN\_OFF} dataset is reduced from 17,450 to 13,016 input-output pairs while the \texttt{BURN\_ON} dataset is reduced from 17,406 to 12,989 pairs.
Each dataset is then rescaled so that the min and max in each coordinate are 0 and 1, respectively, including the output coordinate. 
Rescaling in this way removes bias due to choice of units.


To construct a query point lattice that lies inside the convex hull of the input points, we compute the 25th and 75th percentiles for each of the five input coordinates.
These percentiles define an interval for each dimension, from which we can build a lattice centered around the geometric mean of the inputs.
We use five points per dimension, for a total of $5^5 = 3125$ query points. 
The notion of \qpdf~does not extend directly to these datasets since the inputs were, by design, not drawn uniformly from a bounding box but instead clustered closer to the mean of the inputs.

\renewcommand{\arraystretch}{2}
\begin{figure}[t]
	\centering
	\footnotesize
	\begin{tabular}{ccc}
	\includegraphics[width=.35\textwidth]{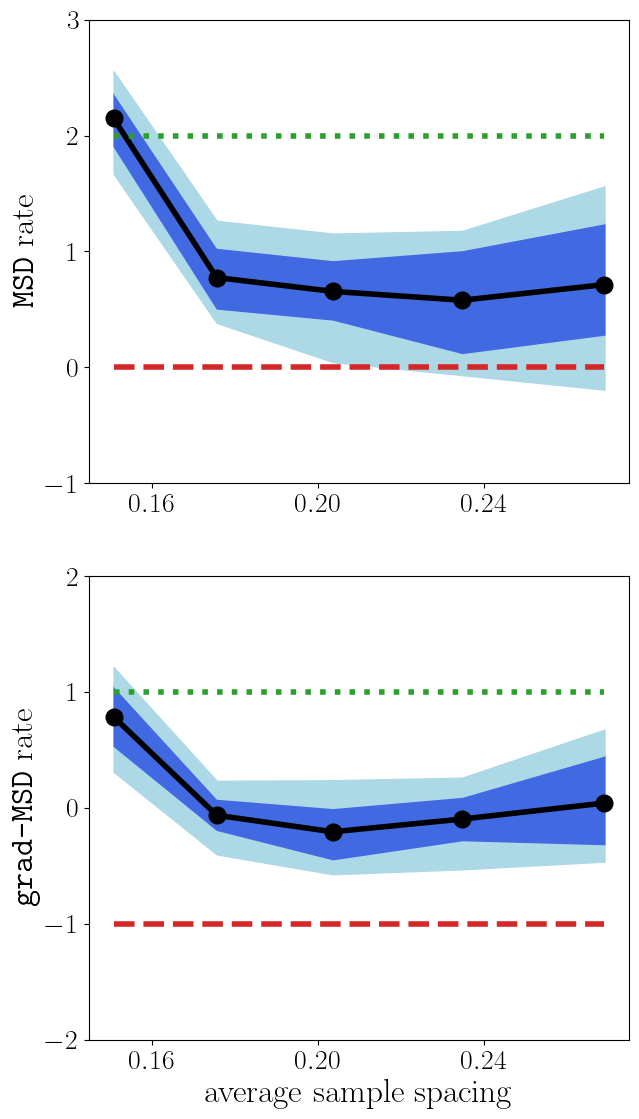}
	\qquad & \qquad
	\includegraphics[width=.35\textwidth]{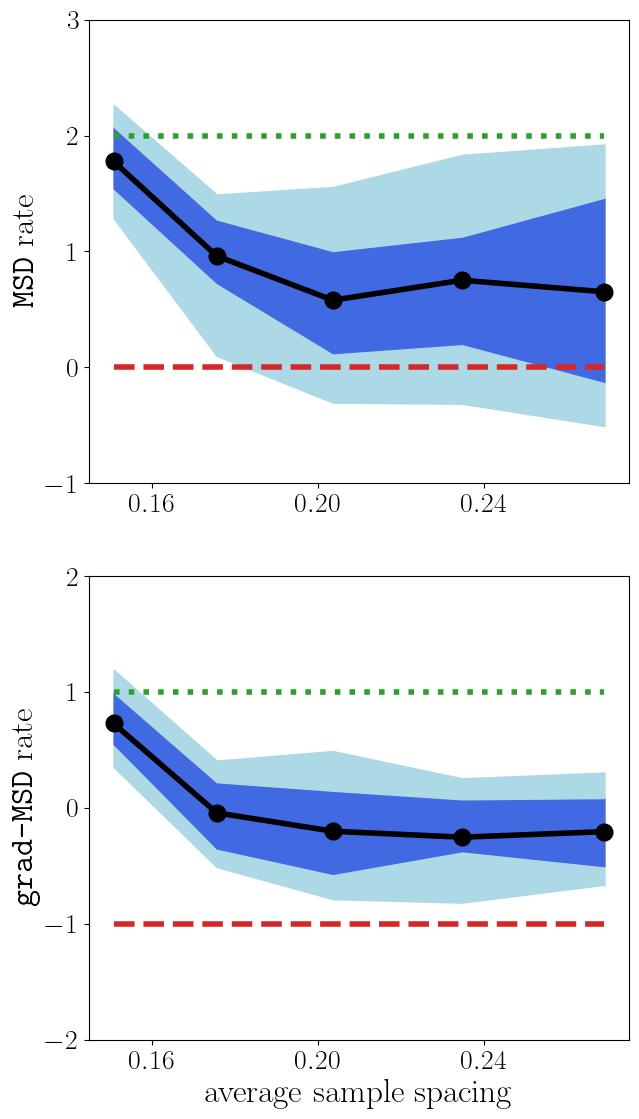}
	\qquad & \qquad
	\includegraphics[width=.15\textwidth]{figs/density_profile_legend.png}
	\\
	 data from \texttt{BURN\_OFF} 
	 \qquad & \qquad
	 data from \texttt{BURN\_ON} 
	\end{tabular}
	\caption{We apply~\cref{alg:DD_static} to two static datasets of $\approx 13,000$ data points each, gathered from simulations of inertial confinement fusion, produced by the HYDRA codebase. The computed rates suggest that the datasets sufficiently capture small scale geometric features, since the rates are nearly equal to the recoverable features rate at the smallest possible \as~value.  In addition, using less than the full data set to create a surrogate model risks introducing spurious features and potentially akin to modeling pure noise. Decreasing \as~further would potentially uncover more features, but would require generating more data, demonstrating how the Delaunay density diagnostic can be used to inform the need for data collection.}
	\label{fig:cogsimdim5}
	\vspace{-.2in}
\end{figure}

In~\cref{fig:cogsimdim5}, we show the results of our experiments for $b=1.2$, $\min n_k = 200$, $\max n_k = 13,016$ for \texttt{BURN\_OFF},$\max n_k = 12,989$ for \texttt{BURN\_ON}, and \# seeds = 50.
With these options, the wall clock time was approximately 2.3 minutes per trial using the same compute resources that produced~\cref{fig:grwk2d3d4d}.
For both datasets, both rates display trends similar to the Griewank function on $[-10,10]^2$, as studied in~\cref{fig:validation1}(b).
Accordingly, the results suggest that the dataset has sufficient density \textit{only} at the smallest possible~\as~value to resolve the features present in the domain of consideration, i.e.~the five dimensional volume of parameters space enclosed by the 25th to 75th percentiles of each coordinate.
The smallest~\as~value occurs when using \textit{all} of the available data, meaning using anything less than the full dataset risks confusing geometric features for noise.

Finally, we note that both datasets are modeling the same physical process: an inertial confinement fusion implosion.
The \texttt{BURN\_OFF} dataset artificially turns off the deposition of energy generated by the fusion reactions, while the \texttt{BURN\_ON} dataset keeps the deposition in place.
The \texttt{BURN\_ON} dataset thus experiences ``yield amplification'' \cite{LHLCB2018,CBL2019} relative to the \texttt{BURN\_OFF} dataset.
Such a process is unlikely to generate features at scales smaller than those present in the \texttt{BURN\_OFF} dataset, though it can change the structure at larger scales. 
Comparing the rate computations between the two datasets in~\ref{fig:cogsimdim5}, we confirm that small scale features are accurately resolved in both, while the \texttt{BURN\_ON} dataset exhibits more variation in rates at lower sample densities (i.e., larger scales).
The differences are not extreme and the similarity in trends between the two datasets serves as a validation of the correlation between the simulations that generated them.

\section{Analytical results} 
\label{sec:analysis}
Our numerical results provide evidence for~\cref{claim:rate2}, \cref{claim:rate0}, \cref{claim:gradrate1}, and \cref{claim:gradrateminus1}.
In this section, we sketch out theoretical support for these claims by employing order-of-magnitude analysis for the quantities appearing in~\cref{alg:DD_msd}. We follow that up with a more formal argument. Throughout the section we assume sample points (those supporting the Delaunay triangulation) and query points (where we evaluate the Delaunay interpolant) are uniformly drawn from the bounding box. Our results are technically valid in the limit of a large number of draws, as we neglect errors coming from finite sampling. In practice, we have found that our results hold as long as deviations from these assumptions are not too significant (e.g. the data does not lie on a lower-dimensional submanifold).

We start with the case of smooth functions, where the sampling resolves all of their features (\cref{claim:rate2}). 
For such functions, the linear interpolant is a reasonable approximation up to quadratic order. 
In particular the difference between successive interpolants will behave as $\|\|\hat{f}_k-\hat{f}_{k-1}\|\|\sim\mathcal{O}(\as_k^2)$, where $\|\|\cdot\|\|$ refers to a discrete approximation to any $L^p$ norm $\left(\text{e.g.,~}\sqrt{\msdi{\bullet}}\right)$ and $\as_k = L/n_k^{1/d}$ is the average sample spacing (the subscript $k$ reminds us of the iteration for which it is computed).
Taking ratios of norms on successive iterations gives us
\begin{equation}\label{eq:smooth_rate_heuristic}
\frac{\|\|\hat{f}_{k-1}-\hat{f}_{k-2}\|\|}{\|\|\hat{f}_{k}-\hat{f}_{k-1}\|\|}\sim\frac{\as_{k-1}^2}{\as_{k}^2}\sim \left(\frac{n_{k}}{n_{k-1}}\right)^{2/d}\sim b^2,
\end{equation} where the last step uses the approximate relation between the number of points and the upsample rate $b$ (i.e.,~line~\ref{algDD:line:upsample} from Algorithm~\ref{alg:DD_msd}). 
Taking $\log_b(\cdot)$ of the expressions in (\ref{eq:smooth_rate_heuristic}) gives $r_k\sim 2$, verifying Claim \ref{claim:rate2}.

With additional notation, we can further formalize the above argument. 
We write the difference between the linear interpolant $\hat{f}_k$ (at iteration $k$) and the true function $f$ as 
\begin{equation}
\hat{f}_k(\bx)-f(\bx) = \sum_{n=2}^\infty a_n^{(k)}(\bx;\sigma_x^{(k)})\langle\Delta x\rangle^n_{\sigma_x^{(k)}}.
\end{equation} 
Here the point $\bx$ is in the $d$-dimensional simplex $\sigma_x^{(k)}$ of the Delaunay triangulation, the discrepancy coefficients $a_n^{(k)}$ depend implicitly on the simplex containing $\bx$, and the average edge length in $\sigma_x^{(k)}$, $\langle\Delta x\rangle_{\sigma_x^{(k)}}$, has been factored out. 
This construction, together with our uniform sampling assumption, guarantees (in a suitable average sense) that the $a_n^{(k)}$ depend only weakly on the refinement iteration $k$, $a_n^{(k)}\approx a_n^{(k-1)}$.

We now write $\langle\Delta x\rangle_{\sigma_x^{(k)}}=w^{(k)}(x)\as_k$, which introduces the weight function $w^{(k)}(x)$ that connects the overall average edge length $\as_k$ in the Delaunay mesh to the local average edge length within $\sigma^{(k)}_x$.
The weight function will be larger in regions of sparse sampling (large simplices) and smaller in regions of dense sampling (small simplices), averaging out to $1$, in a suitable sense. As with the $a_n^{(k)}$, $w^{(k)}(x)$ will depend only weakly on $k$ as a result of our uniform sampling assumption.

We are now ready to construct a discrete $L^1$ distance between successive linear interpolants using the query points $\{\bq_j\}$; a similar argument holds for any discrete $L^p$ norm.
We have:
\begin{align}\label{eq:L1_full}
\|\|\hat{f}_k-\hat{f}_{k-1}\|\| &= \frac{1}{\|\{\bq_j\}\|}\sum_j\|\hat{f}_k(\bq_j)-\hat{f}_{k-1}(\bq_j)\| \nonumber\\
&= \frac{1}{\|\{\bq_j\}\|}\sum_j\left\|\sum_{n=2}^\infty [a_n^{(k)}(\bq_j)w^{(k)}(\bq_j)-b^na_n^{(k-1)}(\bq_j)w^{(k-1)}(\bq_j)]\as_k^n\right\|,
\end{align} where we have used $\as_{k-1}\approx b\as_k$. Since $a_n^{(k)}$ and $w^{(k)}$ (on average) depend only weakly on $k$, we expect the only $k$-dependence of $\|\|\hat{f}_k-\hat{f}_{k-1}\|\|$ to come from $\as_k$. This can be used to more rigorously define the average sense in which $a_n^{(k)}$ and $w^{(k)}$ are independent of $k$ (though we do not attempt that here). If we further assume the $n=2$ term dominates the inner sum, we can conclude that
\begin{equation}\label{eq:smooth_rate_full}
\frac{\|\|\hat{f}_{k-1}-\hat{f}_{k-2}\|\|}{\|\|\hat{f}_{k}-\hat{f}_{k-1}\|\|}\approx\frac{\as_{k-1}^2}{\as_{k}^2}\approx b^2,
\end{equation} confirming our initial estimate from Equation \eqref{eq:smooth_rate_heuristic}.

A similar argument can be made using the gradient of the linear interpolant (and its corresponding norm).
The discrepancy between linear interpolants would scale linearly rather than quadratically.
This causes the ratio of norms to be approximately $b$, and the log-rate to be $1$, as required for \cref{claim:gradrate1}.

We now turn to the case of noisy or highly-oscillatory functions (\cref{claim:rate0}). With insufficient sampling, the linear interpolant will be a poor approximation to the function and will oscillate from iteration to iteration with the amplitude of the noise, which we write as $A_f$:
\begin{equation}
\|\hat{f}_{k}(\bq_j)-\hat{f}_{k-1}(\bq_j)\|\approx A_f\implies \frac{\|\|\hat{f}_{k-1}-\hat{f}_{k-2}\|\|}{\|\|\hat{f}_{k}-\hat{f}_{k-1}\|\|}\approx 1.
\end{equation} Taking the log of the above expression to compute the rate gives us a rate of $0$, verifying \cref{claim:rate0}.

A similar argument holds for the rate of the gradient norm.
In this case, the gradient estimate gets \textit{worse} with refinement, scaling as $A_f/\as_k$.
Thus, the ratio of successive norms will go as $1/b$ and $\log_b$ of the ratio will go to $-1$, as required for~\cref{claim:gradrateminus1}. \\

\paragraph{Analysis of upsampling factor $b$}
Finally, we comment on the trends with upsampling $b$ uncovered in \cref{fig:ackley2}. 
Equation~\eqref{eq:smooth_rate_full} deviates from $b^2$ due to the presence of cubic (and higher order) terms in~\eqref{eq:L1_full}. If we write this deviation as $b^2(1+\varepsilon)$, then the rate becomes
\begin{equation}
\label{eq:rkdev}
r_k~=~\log_b[b^2(1+\varepsilon)]~=~2+\log_b(1+\varepsilon)~\approx~2+\frac{\varepsilon}{\ln b},
\end{equation} where we have made the approximation that the deviation $\varepsilon$ is small. 

From~\eqref{eq:rkdev}, we see that the deviation from the recoverable features \texttt{MSD}-rate (i.e., $2$) gets worse as $b\to 1$. 
This causes stronger fluctuations around the estimated rate for smaller $b$, as seen in \cref{fig:ackley2}.
We can approximate how much additional sampling would be needed to meaningfully reduce these fluctuations. 
Suppose we are upsampling near $1$, so that $b=1+\delta$. 
Since the deviation $\varepsilon$ is sourced by a cubic correction to the rate, we can estimate it as $\varepsilon\approx (a_3/a_2)\as$ (the leading order contribution, $b^2$, already involved two powers of $\as$, leaving just one power remaining in $\varepsilon$). 
If we want small fluctuations around $2$ for the rate, we would need $\varepsilon/\ln b\ll 2$, which simplifies approximately to $\as\ll 2\delta a_3/a_2$.
So, reducing $\delta$ by a factor of two would require a factor of $2^d$ greater total sampling (halving the average sample spacing $\as$) to get the same level of fluctuations around the target rate. 
Thus, our analysis confirms that fluctuations in the computed rates can be controlled by increasing $b$, increasing $n_k$, or both.

This analysis shows that ideally we would select a large value of $b$ to reduce fluctuations, but this choice comes with two drawbacks. First, a large $b$ makes it difficult to resolve large scales. The diagnostic requires two iterations to compute the first rate, so if we start with an initial spacing of $\as_0=L/n_0^{(1/d)}$, the first diagnostic point will be at a spacing of $\as_2=\as_0/b^2$, which might be significantly smaller than the box side $L$. Unfortunately, there is little that can be done to alleviate this loss of large scales. Fluctuations only decrease with increasing $b$ or decreasing spacing, both of which lose track of large scales. Secondly, a large $b$ means fewer points on the diagnostic curve. This will give us a less precise estimate of how much data is too much (if the diagnostic shows features can be resolved). This is less of an issue if we primarily care about the ability of the full dataset to resolve features. In that case, the only concern with large $b$ is having a sufficiently large initial sample $n_0$ to avoid large fluctuations. This depends somewhat on the structure of the function and is best determined empirically.

\section{Conclusions, extensions, and code} 
We have demonstrated in this paper how the convergence rate of iteratively refined, piecewise linear approximations to a function can be used to assess if the function has been sampled densely enough \textit{relative} to its variation.
Our computational technique eschews nearly any assumption on $f$ as it detects pure noise and undersampled oscillations as equivalent phenomena.
Many extensions of the approach are plausible, including assessing convergence in other norms (as mentioned previously) and consideration of higher dimensional data sets, additional static data sets, discontinuous functions, functions with singularities, time-dependent functions, and so forth.

To aid any interested parties in exploring these and other directions, we have released Python code and two driver scripts that effectively replicate the numerical results shown in~\cref{fig:grwk2d3d4d} for $f=g_2$ and in \cref{fig:topography} for the static topography dataset.
The parameters were adjusted so that the requisite data and figures can be generated in minutes using a typical laptop with a standard modern python environment.
The code is included as part of the published version of this work.
A version of the code is also available on Github~\cite{dddrepo}.



\section*{Acknowledgments}
We would like to thank Peer-Timo Bremer, Luc Peterson, Brian Spears, and many other colleagues at Lawrence Livermore National Laboratory for providing feedback, encouragement, and financial support for this research.
We would also like to thank Tyler Chang at Argonne National Laboratory for many helpful discussions regarding the \texttt{DelaunaySparse} software.
Finally, we would like to thank the anonymous reviewers of the manuscript whose suggestions greatly improved the manuscript and accompanying code.
This work was performed under the auspices of the U.S.~Department of Energy by Lawrence Livermore National Laboratory under Contract DE--AC52--07NA27344 and the LLNL-LDRD Program under Project tracking No.\ 21--ERD--028.  
Release number LLNL--JRNL--846568. 


\bibliographystyle{ACM-Reference-Format}
\bibliography{main.bib}

\appendix
\section{Computing the gradient of the Delaunay interpolant}
\label{app:grad}
We can compute $\nabla \hat f_k(\bq_i)$ in~\cref{alg:DD_grad-msd} by exploiting some standard techniques of differential geometry and linear algebra.
Since $\hat f_k:CH(\calD)\raw\bbR^1$, and $CH(\calD)\subset\bbR^d$, we can write
\[
x_{d+1} = \hat f_k(x_1,\ldots,x_d),
\]
where $\{x_i\}$ are standard Euclidean coordinates for $\bbR^{d+1}$.  
Thus,
\begin{equation}
\label{eq:impsfc}
\hat f_k(x_1,\ldots,x_d) - x_{d+1} = 0
\end{equation}
defines an implicit $d$-dimensional surface in $\bbR^{d+1}$ that is piecewise flat.
Define $\hat n_k(\bq_i)$ by implicit differentiation as
\begin{equation}
\label{eq:n_to_grad}
\hat n_k(\bq_i) = \left( \partial_1 \hat f_k(\bq_i), \cdots, \partial_d \hat f_k(\bq_i), -1\right) = \left(\nabla \hat f_k(\bq_i), -1\right).
\end{equation}
Observe $\hat n_k(\bq_i)/\|\|\hat n_k(\bq_i)\|\|$ is a unit normal vector to the piecewise flat implicit surface (\ref{eq:impsfc}).

Now, near $\bq_i$, the surface (\ref{eq:impsfc}) is determined by the values of the Delaunay $d$-simplex $\calS$ that contains $\bq_i$.
Let $\{\bs_1,\ldots,\bs_{d+1}\}$ denote the vertices of $\calS$, which are found during the computation of ${\hat f}_k(\bq_i)$ by~\cref{alg:DS}.
Recall that each $\bs_\ell\in\{\bx_j\}$, meaning $f(\bs_\ell)$ is known by assumption.
Thus, $\{(\bs_\ell, f(\bs_\ell))\}$ is a collection of $d+1$ points in $\bbR^{d+1}$ lying on the surface (\ref{eq:impsfc}) near $\bq_i$.

Set $A$ to be the $(d+1)\times (d+1)$ matrix whose rows are formed by the vectors $(\bs_\ell, f(\bs_\ell))$.
Subtract the column-wise average of $A$ from each row of $A$, which has the effect of translating the barycenter of $\calS$ to the origin.
Let $U\Sigma V^\ast$ be the SVD of $A$.
Then $V^\ast$ is an orthonormal set whose first $d$ vectors form a basis for (\ref{eq:impsfc}), meaning the last vector, call it $\bv_{d+1}$, is a unit normal to (\ref{eq:impsfc}).
By scaling $\bv_{d+1}$ so that its last coordinate is $-1$, we have found $\hat n_k(\bq_i)$ and can recover $\nabla \hat f_k(\bq_i)$ from (\ref{eq:n_to_grad}).

We remark briefly on the case where $\bq_i$ lies at the interface of one or more Delaunay simplices.
The gradient is not continuous across mesh elements (since the function is piecewise flat)  and thus $\nabla\hat f_k(\bq_i)$ has no unique definition. 
In practice, however, this is very unlikely to occur and such instances could be detected and managed robustly in a number of ways.







\end{document}